\documentclass[pdflatex,sn-mathphys-num]{sn-jnl}

\usepackage[justification=raggedright,singlelinecheck=true,tableposition=top]{caption} 
\usepackage{graphicx}
\usepackage{subfigure}
\usepackage{epstopdf}
\usepackage{mathtools}

\usepackage{amssymb}

\usepackage{amsmath}

\usepackage{amsfonts}
\usepackage{amsthm}
\usepackage{mathrsfs}

\usepackage{marginnote}
\usepackage{color,soul}
\usepackage{setspace}
\usepackage{multirow}
\usepackage{url}
\usepackage{colortbl} 
\usepackage{bm} 
\usepackage{enumerate}
\usepackage{isomath} 
\usepackage{setspace}
\usepackage{dblfloatfix} 

\usepackage[title]{appendix}%
\usepackage{xcolor}%
\usepackage{textcomp}%
\usepackage{manyfoot}%
\usepackage{booktabs}%
\usepackage{algorithm}%
\usepackage{algorithmicx}%
\usepackage{algpseudocode}%
\usepackage{listings}%

\newcommand{\eref}[1]{(\ref{#1})}
\newcommand{\fref}[1]{Fig.~\ref{#1}}

\newcommand{\vm}[1]{\bm{#1}} 
\newcommand{\bsym}[1]{\bm{#1}}

\renewcommand{\Re}{{\rm{I\!R}}}
\newcommand{\diffx}{d\vm{x}}
\newcommand{\diffs}{ds}
\newcommand{\transpose}{\mathsf{T}}

\newcommand{\cons}{\mathrm{c}}
\newcommand{\stab}{\mathrm{s}}
\newcommand{\trace}[1]{\mathrm{trace}{\,\left(#1\right)}}
\newcommand{\smat}[2][ccccccccccccccccccccccccccccccccccccccccccccccccccc]{\left
[\arraycolsep=3.0pt\def\arraystretch{1.2}\begin{array}{#1}#2 \\ \end{array} \right]}

\newcommand{\cmat}[2][ccccccccccccccccccccccccccccccccccccccccccccccccccc]{\left
\{\arraycolsep=0.4pt\def\arraystretch{1.2}\begin{array}{#1}#2 \\ \end{array} \right\}}

\newcommand{\alejandro}[1]{{\color{black}{#1}}}

\theoremstyle{thmstyleone}%
\theoremstyle{thmstyletwo}%

\theoremstyle{thmstylethree}%

\begin{document}

\title[A node-based uniform strain virtual element method for elastoplastic solids]{A node-based uniform strain virtual element method for elastoplastic solids}


\author[1,2]{\fnm{Rodrigo} \sur{Silva-Valenzuela}}\email{rosilva@ug.uchile.cl}

\author*[1]{\fnm{Alejandro} \sur{Ortiz-Bernardin}}\email{aortizb@uchile.cl}

\author[3]{\fnm{Edoardo} \sur{Artioli}}\email{eartioli@unimore.it}


\affil[1]{\orgdiv{Computational and Applied Mechanics Laboratory, Department of Mechanical Engineering}, 
           \orgname{Universidad de Chile}, 
           \orgaddress{\street{Av. Beauchef 851}, \city{Santiago}, \postcode{8370456}, \country{Chile}}}

\affil[2]{\orgdiv{Department of Mechanical Engineering}, 
          \orgname{Universidad de La Serena}, 
          \orgaddress{\street{Av. Benavente 980}, \city{La Serena}, \postcode{1720170}, \country{Chile}}}

\affil[3]{\orgdiv{Dipartimento di Ingegneria ``Enzo Ferrari''}, 
          \orgname{Universit\`a degli Studi di Modena e Reggio Emilia}, 
          \orgaddress{\street{Via Pietro Vivarelli 10}, \city{Modena}, \postcode{41125}, \country{Italy}}}



\abstract{A recently proposed node-based uniform strain virtual element method (NVEM)
is here extended to small strain elastoplastic solids. In the proposed method, 
the strain is averaged at the nodes from the strain of surrounding 
linearly precise virtual elements using a generalization to virtual elements 
of the node-based uniform strain approach for finite elements. 
The averaged strain is then used to sample the weak form at the nodes 
of the mesh leading to a method in which all the field variables, 
including state and history-dependent variables, 
are related to the nodes and thus they are tracked only at these locations
during the nonlinear computations. Through various elastoplastic
benchmark problems, we demonstrate that the NVEM is locking-free 
while enabling linearly precise virtual elements to solve 
elastoplastic solids with accuracy.}

\keywords{Virtual element method, Nodal integration, Strain averaging, Uniform strain, Volumetric locking, Elastoplasticity}



\maketitle

\section{Introduction}
\label{sec:intro}

Elastoplastic solids demand advanced numerical techniques when using
Galerkin-based approaches such as the finite element method (FEM)
and the virtual element method (VEM). The need for advanced techniques
has its roots in the presence of volumetric locking in the numerical solution
due to the volume preserving nature of the plastic strain (plastic 
incompressibility condition) and in the volume preserving condition 
that arises when the Poisson's ratio approaches $1/2$ (elastic 
incompressibility condition). Whichever is the source of the locking 
behavior, the standard lowest order elements in general perform poorly 
in problems that involve volume preserving conditions. 
In FEM, several approaches to deal with locking effects are found
in the literature. An exhaustive review of these approaches is out of the scope
of this paper, but we mention the most relevant ones: 
reduced/selective integration~\cite{Malkus-Hughes:1978}, 
B-bar technique~\cite{Hughes:1980,Simo-Hughes:1998}, 
mixed formulations~\cite{Malkus-Hughes:1978}, 
assumed strain methods~\cite{Simo-Rifai:1990}, and
nodal integration~\cite{Bonet-Burton:1998,Dohrmann-Heinstein-Jung-Key-Witkowski:2000,
Bonet-Marriott-Hassan:2001,Puso-Solberg:2006,Puso-Chen-Zywicz-Elmer:2008,Krysl-Zhu:2008,
Broccardo-Micheloni-Krysl:2009,Castellazzi-Krysl:2012,Krysl-Kagey:2012,Artioli-Castellazzi-Krysl:2014}. 
Of particular interest for the method proposed
in this paper are nodal integration techniques. In these approaches,
the Galerkin weak form is sampled at the nodes of the mesh
leading to methods in which all the field variables (including state and history-dependent
variables) are associated with the nodes. 

The VEM is a generalization of the FEM to elements with arbitrary number of edges/faces 
(convex or nonconvex polytopes) known as virtual elements~\cite{BeiraoDaVeiga-Brezzi-Cangiani-Manzini-Marini-Russo:2013}.
In its standard form, the method consists in the construction of an algebraic
(exact) representation of the stiffness matrix without computation of basis
functions (basis functions are \textit{virtual}). In this process, a decomposition
of the stiffness matrix into a consistency part and a stability part that
ensures convergence of the method~\cite{Cangiani-Manzini-Russo-Sukumar:2015}
is realized. The VEM has gained much interest in recent years 
and nowadays its applications can be found, for instance, in
elastic and inelastic solids~\cite{BeiraodaVeiga-Brezzi-Marini:2013,
BeiraodaVeiga-Lovadina-Mora:2015,Artioli-BeiraoDaVeiga-Lovadina-Sacco:2017,
Artioli-BeiraoDaVeiga-Lovadina-Sacco:2017b,Park-Chi-Paulino:2021,
Gain-Talischi-Paulino:2014,Daltri-deMiranda-Patruno-Sacco:2021,
Tang-Liu-Zhang-Feng:2020,Cihan-Hudobivnik-Aldakheel-Wriggers:2021}, 
elastodynamics~\cite{Park-Chi-Paulino:2020,Park-Chi-Paulino:2021b}, 
finite deformations~\cite{Wriggers-Rust:2019,DeBellis-Wriggers-Hudobivnik:2019,
Aldakheel-Hudobivnik-Wriggers:2019,Wriggers-Hudobivnik:2017,Wriggers-Reddy-Rust-Hudobivnik:2017,
Hudobivnik-Aldakheel-Wriggers:2019,Zhang-Chi-Paulino:2020,Chi-BeiraoDaVeiga-Paulino:2017,
vanHuyssteen-Reddy:2020}, contact mechanics~\cite{Wriggers-Rust:2019,Wriggers-Rust-Reddy:2016,
Aldakheel-Hudobivnik-Artioli-BeiraoDaVeiga-Wriggers:2020}, 
fracture mechanics~\cite{Aldakheel-Hudobivnik-Wriggers:2019b,
Hussein-Aldakheel-Hudobivnik-Wriggers-Guidault-Allix:2019,
NguyenThanh-Zhuang-NguyenXuan-Rabczuk-Wriggers:2018,
Benedetto-Caggiano-Etse:2018,Artioli-Marfia-Sacco:2020}, 
fluid mechanics~\cite{BeiraoDaVeiga-Pichler-Vacca:2021,Chen-Wang:2019,
Gatica-Munar-Sequeira:2018,BeiraoDaVeiga-Lovadina-Vacca:2018,
Chernov-Marcati-Mascotto:2021}, geomechanics~\cite{Andersen-Nilsen-Raynaud:2017,
Lin-Zheng-Jiang-Li-Sun:2020} and 
topology optimization~\cite{Zhang-Chi-Paulino:2020,Paulino-Gain:2014,
Chi-Pereira-Menezes-Paulino:2020}. 

In the VEM literature, there are few methods already developed that are suitable
for modeling nearly incompressible elastic solids. These are generalizations
of some of the above mentioned approaches for FEM. For instance,
B-bar formulation~\cite{Park-Chi-Paulino:2020,Park-Chi-Paulino:2021},
mixed formulation~\cite{Artioli-deMiranda-Lovadina-Patruno:2017}, enhanced strain
formulation~\cite{Daltri-deMiranda-Patruno-Sacco:2021}, hybrid formulation~\cite{Dassi-Lovadina-Visinoni:2021},
and nonconforming formulations~\cite{Kwak-Park:2022,Zhang-Zhao-Yang-Chen:2019,Yu:2023}.
In small strain elastoplasticity, the VEM literature is very limited. For instance,
in Refs.~\cite{Artioli-BeiraoDaVeiga-Lovadina-Sacco:2017b,Xu-Wang-Wriggers:2024} 
high-order VEM has been used to improve
the numerical performance when facing locking effects in elastoplastic 
solids. A mixed formulation based on 
a Hu-Washizu functional was adopted in Ref.~\cite{Cihan-Hudobivnik-Aldakheel-Wriggers:2021}. 
Recently, a stabilization-free hybrid virtual element method has
been proposed for elastoplastic solids that is locking-free~\cite{Liguori-Madeo-Marfia-Garcea-Sacco:2024}.

In this paper, the recently proposed node-based uniform strain 
virtual element method~\cite{Ortiz-Silva-Salinas-Hitschfeld-Luza-Rebolledo:2023} (NVEM) 
is extended to small strain elastoplastic solids. In the proposed approach, 
the strain is averaged at the nodes from the strain of surrounding 
linearly precise virtual elements using a generalization to 
virtual elements of the node-based uniform strain approach for 
finite elements~\cite{Dohrmann-Heinstein-Jung-Key-Witkowski:2000}.
The nodal strain that results from the averaging process is interpreted 
as the nodal sample of the strain in the nodal integration of the
weak form. Consequently, the nodal strain is also used at the constitutive
evaluation level. As in any nodal integration method, the state 
and history-dependent variables in the NVEM become associated with the nodes.
In practice, this means that in nonlinear computations these variables
are tracked only at the nodes. This feature can be exploited to avoid mesh 
remapping of these variables in Lagrangian large deformation simulations with remeshing 
(see, for instance, Ref.~\cite{Puso-Chen-Zywicz-Elmer:2008}), which is not possible
when Gauss integration is used. We do not intend to explore
the latter feature in this paper as the focus here is on the small
strain regime, which is a necessary intermediate step towards developments
in the finite strain regime.

The remainder of this paper is structured as follows. In Section~\ref{sec:nvem},
the NVEM for small strain elastoplasticity is developed. The elastoplastic
constitutive model used and the stabilization proposed for the NVEM are
presented in Section~\ref{sec:stabilization}. In Section~\ref{sec:numexp}, 
various elastoplastic benchmark problems are considered to assess the performance
of the NVEM. The paper ends with a summary and conclusions in Section~\ref{sec:conclusions}.

\section{Node-based uniform strain virtual element method}
\label{sec:nvem}

In this section, the basics of the node-based uniform strain virtual element
method (NVEM), which was developed in 
Ref.~\cite{Ortiz-Silva-Salinas-Hitschfeld-Luza-Rebolledo:2023}, is summarized. 
The method belongs to the Galerkin weak formulation family of methods. 
In this sense, we consider an elastic body that occupies the
open domain $\Omega \subset \Re^2$ and is
bounded by the one-dimensional surface $\Gamma$ whose
unit outward normal is $\vm{n}_{\Gamma}$. The boundary is assumed
to admit decompositions $\Gamma=\Gamma_D\cup\Gamma_N$ and
$\emptyset=\Gamma_D\cap\Gamma_N$, where $\Gamma_D$ is the Dirichlet
boundary and $\Gamma_N$ is the Neumann boundary. The closure of
the domain is $\overline{\Omega}=\Omega\cup\Gamma$. Let
$\vm{u}(\vm{x}) : \overline{\Omega} \rightarrow \Re^2$ be
the displacement field at a point of the elastic body 
with position vector $\vm{x}$ when the body is subjected to external tractions
$\vm{t}_N(\vm{x}):\Gamma_N\rightarrow \Re^2$ and body forces $\vm{b}(\vm{x}):\Omega\rightarrow\Re^2$.
The imposed Dirichlet (essential) boundary conditions are
$\vm{u}_D(\vm{x}):\Gamma_D\rightarrow \Re^2$. The displacement field
$\vm{u}(\vm{x})\in \mathcal{V}$ is found such that (weak form)
\begin{equation}\label{eq:weakform}
\begin{split}
a(\vm{u},\vm{v}) &= \ell(\vm{v}) \quad \forall \vm{v}(\vm{x})\in \mathcal{W},\\
a(\vm{u},\vm{v}) = \int_{\Omega}\bsym{\sigma}(\vm{u}):\bsym{\varepsilon}(\vm{v})\,\diffx,
& \quad \ell(\vm{v}) = \int_{\Omega}\vm{b}\cdot\vm{v}\,\diffx + \int_{\Gamma_N}\vm{t}_N\cdot\vm{v}\,\diffs,
\end{split}
\end{equation}
where $\mathcal{V}$ denotes the space of admissible displacements 
and $\mathcal{W}$ the space of its variations; $\bsym{\sigma}$ is the Cauchy stress tensor 
and $\bsym{\varepsilon}$ is the small strain tensor that is given by
\begin{equation}\label{eq:symstrain}
\bsym{\varepsilon}(\vm{u})= 
\frac{1}{2}\left(\vm{u}\otimes\bsym{\nabla}+\bsym{\nabla}\otimes\vm{u}\right).
\end{equation}

\subsection{Virtual element method}
\label{sec:vem}

\begin{figure}[!bth]
\centering
\includegraphics[width=0.5\linewidth]{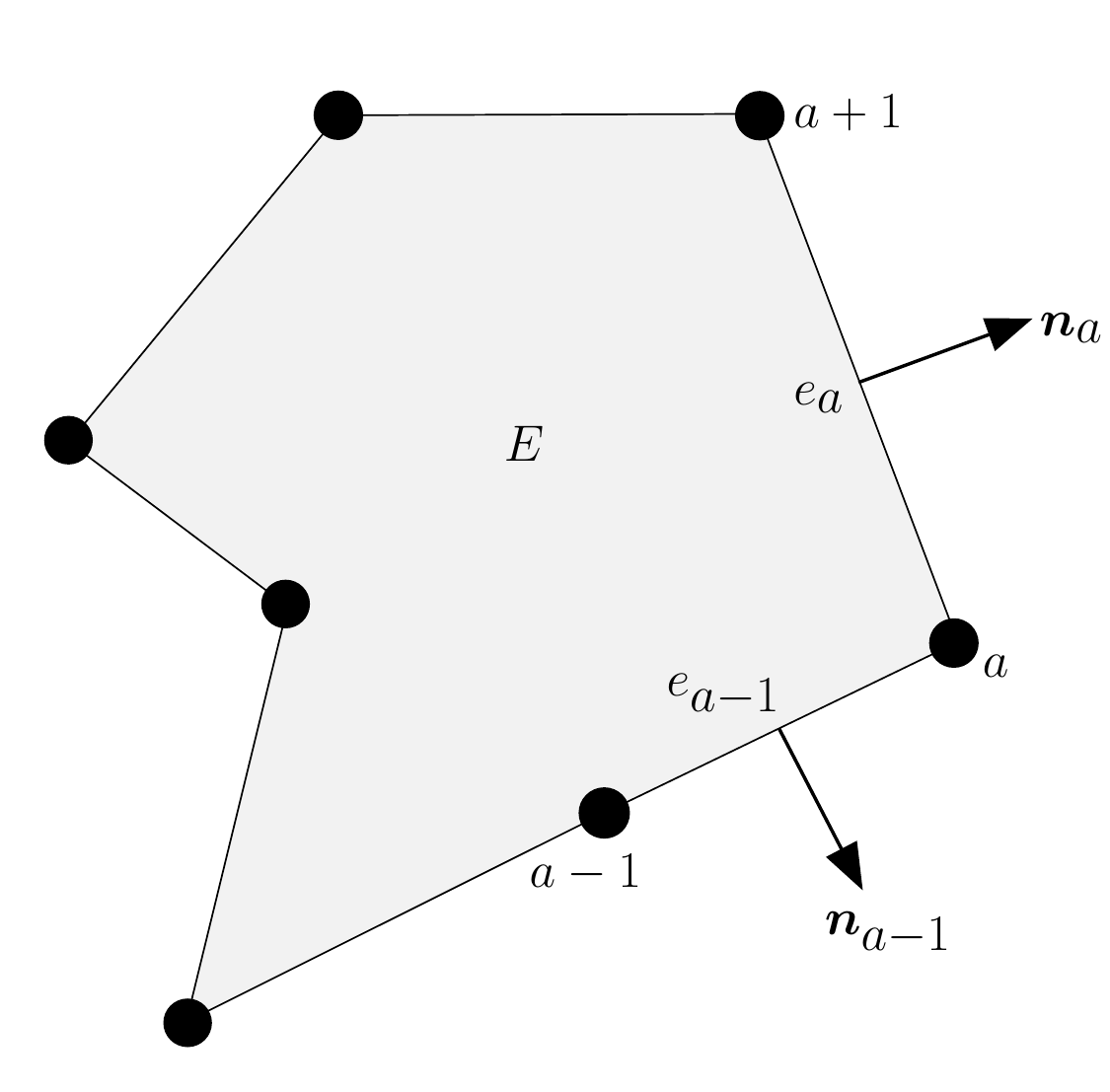}
\caption{Schematic representation of a polygonal element of $N_E^V=7$ edges}
\label{fig:polyelement}
\end{figure}

The weak form~\eref{eq:weakform} is the continuous problem. The discrete problem
is formulated on a partition of the domain $\Omega$ into nonoverlapping elements
with arbitrary number of edges (convex or non-convex polygons). 
This partition is denoted by $\mathcal{T}_h$, where $h$ is the maximum diameter 
of any element in the partition. An element in the partition is denoted by $E$ 
and its boundary by $\partial E$. $|E|$ is the area of the element and $N_E^V$
its number of edges/nodes. The unit outward normal to the element boundary 
in the Cartesian coordinate system is denoted by $\vm{n}=[n_1 \quad n_2]^\transpose$. 
\fref{fig:polyelement} depicts an element with seven edges ($N_E^V=7$), where 
the edge $e_a$ of length $|e_a|$ and the edge $e_{a-1}$ of length $|e_{a-1}|$ are the 
element edges incident to node $a$, and $\vm{n}_a$ and $\vm{n}_{a-1}$ are the unit 
outward normals to these edges, respectively.

Following a standard Galerkin approach, we assume approximations of
$\vm{u}$ and $\vm{v}$ on the element, as follows:
\begin{equation}\label{eq:discrete_field}
\begin{split}
\vm{u}_h &=\cmat{u_{1h} \\ u_{2h}}=\sum_{a=1}^{N_E^V}\phi_a(\vm{x})\vm{u}_a, \quad \vm{u}_a = \cmat{u_{1a} \\ u_{2a}},\\
\vm{v}_h &=\cmat{v_{1h} \\ v_{2h}}=\sum_{a=1}^{N_E^V}\phi_a(\vm{x})\vm{v}_a, \quad \vm{v}_a = \cmat{v_{1a} \\ v_{2a}}.
\end{split}
\end{equation}
where $\{\phi_a(\vm{x})\}_{a=1}^{N_E^V}$ are basis functions that form a partition of unity.
A peculiarity of the VEM is that the basis functions are never computed,
which is why they are considered \textit{virtual}. For the method to work, we only need to
assume their behavior on the element boundary. For linearly precise VEM approximations,
the basis functions on the element boundary are assumed to be
\begin{itemize}
        \item piecewise linear (edge by edge),
        \item continuous on the element edges,
\end{itemize}
which means that the basis functions possess the Kronecker-delta property on the element edges,
and hence they behave like the one-dimensional hat function.

At the element level, the following discrete local spaces are defined:
\begin{equation*}
\mathscr{V}_h|_E :=
\left\{\vm{u}_h(\vm{x}): \vm{u}_h \in \mathcal{V}(E)\right\},\quad
\mathscr{W}_h|_E := \mathscr{V}_h|_E.
\end{equation*}
The discrete local spaces are assembled to form the following discrete global
spaces:
\begin{equation*}
\begin{split}
\mathscr{V}_h &:=
\left\{\vm{u}(\vm{x})\in\mathscr{V}: \vm{u}|_E \in \mathscr{V}_h|_E 
\quad\forall E \in \mathcal{T}_h\right\},\\
\mathscr{W}_h &:=
\left\{\vm{v}(\vm{x})\in\mathscr{W}: \vm{v}|_E \in \mathscr{V}_h|_E 
\quad \forall E \in \mathcal{T}_h\right\}.
\end{split}
\end{equation*}
Using the preceding definitions, the discrete version of the weak form~\eref{eq:weakform} 
reads: find $\vm{u}_h \in\mathscr{V}_h$ such that
\begin{equation}\label{eq:weakform_disc}
\sum_{E\in\mathcal{T}_h} a_E(\vm{u}_h,\vm{v}_h) = 
\sum_{E\in\mathcal{T}_h} \ell_E(\vm{v}_h)\quad\forall\vm{v}_h\in\mathscr{W}_h.
\end{equation}

To obtain the discrete weak form for the VEM,
we follow the standard VEM literature (see for instance, Ref.~\cite{BeiraodaVeiga-Brezzi-Marini-Russo:2014}).
We first define a projection operator $\Pi$ onto the space of polynomials of degree 1. 
To this end, let $[\mathcal{P}(E)]^2$ represent the space of polynomials of degree 1 over 
the element $E$. The  projection operator $\Pi$ is defined as:
\begin{equation}\label{eq:projection_operator}
\Pi: \mathscr{V}_h|_E  \to [\mathcal{P}(E)]^2, \,
\Pi\vm{p}=\vm{p} \quad \forall \vm{p}\in [\mathcal{P}(E)]^2.
\end{equation}
$\Pi$ is then used to split the displacement approximation
on the element, as follows:
\begin{equation}\label{eq:vem_decomp}
\vm{u}_h=\Pi\vm{u}_h+(\vm{u}_h-\Pi\vm{u}_h),
\end{equation}
where $\Pi\vm{u}_h$ is the polynomial part of $\vm{u}_h$ (of degree 1) 
and $\vm{u}_h-\Pi\vm{u}_h$ contains its remainder terms. The remainder 
terms can contain polynomials of order greater than 1 or even nonpolynomial terms.
The actual form of the projection $\Pi\vm{u}_h$ is obtained from the 
orthogonality condition: $\forall\vm{p}\in [\mathcal{P}(E)]^2$,
\begin{equation}\label{eq:ortho_condition}
a_E(\vm{u}_h-\Pi\vm{u}_h,\vm{p}) = a_E(\vm{p},\vm{v}_h-\Pi\vm{v}_h) = 0, 
\end{equation}
which, at the element level, gives~\cite{Gain-Talischi-Paulino:2014,Ortiz-Russo-Sukumar:2017,Ortiz-Alvarez-Hitschfeld-Russo-Silva-Olate:2019,Silva-Ortiz-Sukumar-Artioli-Hitschfeld:2020}
\begin{equation}\label{eq:operator}
\Pi\vm{u}_h =
\smat{(x_1-\bar{x}_1) & 0 & \frac{(x_2-\bar{x}_2)}{2} & 1 & 0 & \frac{(x_2-\bar{x}_2)}{2}\\
       0 & (x_2-\bar{x}_2) & \frac{(x_1-\bar{x}_1)}{2} & 0 & 1 & \frac{(\bar{x}_1-x_1)}{2}}
\cmat{\widehat{\varepsilon}_{11} \\ \widehat{\varepsilon}_{22} 
       \\ 2\,\widehat{\varepsilon}_{12} \\ \bar{u}_{1} \\ \bar{u}_{2}
       \\ 2\,\widehat{\omega}_{12} },
\end{equation}
where $\bar{x}_1$ and $\bar{x}_2$ are the components of the mean of the values that
the position vector $\vm{x}=\left[x_1 \quad x_2\right]^\transpose$ takes over 
the vertices of the element; i.e.,
\begin{equation}\label{eq:bar_x}
\bar{\vm{x}}=\cmat{\bar{x}_1 \\ \bar{x}_2}=\frac{1}{N_E^V}\sum_{a=1}^{N_E^V}\vm{x(\vm{x}_a)},
\end{equation}
where $\vm{x}_a=\left[x_{1a} \quad x_{2a}\right]^\transpose$
are the coordinates of node $a$; $\bar{u}_{1}$ and $\bar{u}_{2}$ are the components 
of the mean of the values that the displacement approximation $\vm{u}_h=\left[u_{1h} \quad u_{2h}\right]^\transpose$ 
takes over the vertices of the element; i.e.,
\begin{equation}\label{eq:bar_disp}
\bar{\vm{u}}=\cmat{\bar{u}_{1} \\ \bar{u}_{2}}=\frac{1}{N_E^V}\sum_{a=1}^{N_E^V}\vm{u}_h(\vm{x}_a),
\end{equation}
In other words, $\bar{\vm{x}}$ and $\bar{\vm{u}}$ represent the geometric center of
the element and its associated displacement vector, respectively;
the terms $\widehat{\varepsilon}_{ij}$ are components of the element average
\begin{equation}\label{eq:avg_strain}
\widehat{\bsym{\varepsilon}}(\vm{u}_h) = \frac{1}{|E|}\int_E \bsym{\varepsilon}(\vm{u}_h)\,\diffx
= \frac{1}{2|E|}\int_{\partial E}\left(\vm{u}_h\otimes\vm{n}+\vm{n}\otimes\vm{u}_h\right)\,\diffs,
\end{equation}
and $\widehat{\omega}_{12}$ is the component of the element average
\begin{equation}\label{eq:avg_skewstrain}
\widehat{\bsym{\omega}}(\vm{u}_h) = \frac{1}{|E|}\int_E \bsym{\omega}(\vm{u}_h)\,\diffx\\
=\frac{1}{2|E|}\int_{\partial E}\left(\vm{u}_h\otimes\vm{n}-\vm{n}\otimes\vm{u}_h\right)\,\diffs,
\end{equation}
where $\bsym{\omega}(\vm{u}_h)$ is the skew-symmetric tensor that represents rotations. 

On substituting~\eref{eq:vem_decomp} into~\eref{eq:weakform_disc}, and using 
the orthogonality condition~\eref{eq:ortho_condition} and noting that
$\vm{u}_h$ and $\vm{v}_h \in [\mathcal{P}(E)]^2$, leads to
the following VEM representation of the discrete weak form: 
find $\vm{u}_h \in\mathscr{V}_h$ such that
\begin{equation}\label{eq:vemweakform_disc}
\sum_{E\in\mathcal{T}_h} \biggl[a_E(\Pi\vm{u}_h,\Pi\vm{v}_h) + s_E(\vm{u}_h-\Pi\vm{u}_h,\vm{v}_h-\Pi\vm{v}_h) \biggr] 
 = \sum_{E\in\mathcal{T}_h} \ell_E(\Pi\vm{v}_h)\quad\forall\vm{v}_h\in\mathscr{W}_h,
\end{equation}
where $s_E(\vm{u}_h-\Pi\vm{u}_h,\vm{v}_h-\Pi\vm{v}_h)$ is a computable approximation to
$a_E(\vm{u}_h-\Pi\vm{u}_h,\vm{v}_h-\Pi\vm{v}_h)$ and is meant to provide stability.

\subsection{Nodal averaging operator}
\label{sec:nodavgoperator}

The VEM as described above is prone to volumetric locking in the 
limit $\nu \to 1/2$. Using the virtual element mesh and considering a typical
nodal vertex $I$, the NVEM that
was proposed in Ref.~\cite{Ortiz-Silva-Salinas-Hitschfeld-Luza-Rebolledo:2023} applies 
a nodal averaging operator $\pi_I$ to~\eref{eq:vemweakform_disc} 
that precludes volumetric locking without introducing 
additional degrees of freedom\footnote{\alejandro{The nodal averaging operator
permits to integrate the weak form integrals directly at the nodes. This
results in a total number of incompressiblity constraints equal to 
the number of nodes in the mesh. If this number divides the total 
number of displacement equations, two degrees of freedom every 
one constraint is obtained in two dimensions, which is the optimal ratio to 
perform well in incompressible and nearly incompressible 
settings~\cite{hughes:2000}.}}. This leads to a nodal version of~\eref{eq:vemweakform_disc},
as follows:
find $\vm{u}_h \in\mathscr{V}_h$ such that
\begin{equation}\label{eq:nodal_vemweakform_disc}
\begin{split}
\sum_{I\in\mathcal{T}_h} &\biggl[a_I\bigl(\pi_I[\Pi\vm{u}_h],\pi_I[\Pi\vm{v}_h]\bigr)
+ s_I\bigl(\pi_I[\vm{u}_h-\Pi\vm{u}_h],\pi_I[\vm{v}_h-\Pi\vm{v}_h]\bigr) \biggr]\\
&= \sum_{I\in\mathcal{T}_h} \ell_I\bigl(\pi_I[\Pi\vm{v}_h]\bigr)
\quad\forall\vm{v}_h\in\mathscr{W}_h,
\end{split}
\end{equation}
where the notations $a_I$, $s_I$ and $\ell_I$ are
introduced as the nodal counterparts of $a_E$, $s_E$ and $\ell_E$, respectively.
The construction of the nodal averaging operator is described next.

\begin{figure}[!bth]
\centering
\includegraphics[width=0.5\linewidth]{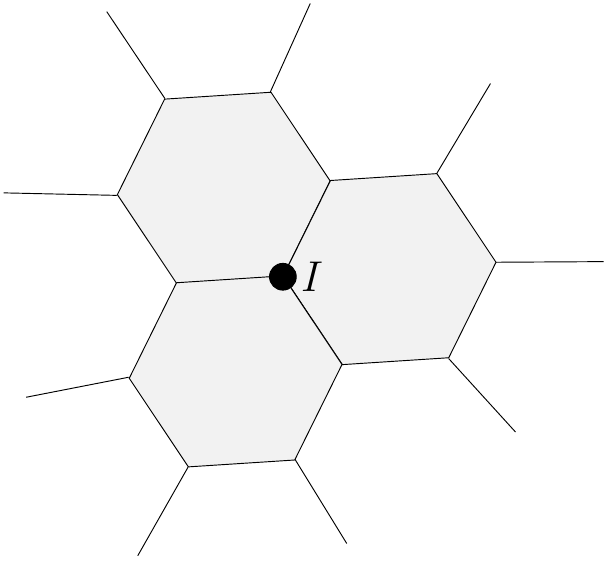}
\caption{Nodal patch $\mathcal{T}_I$ (shaded elements) formed by the virtual elements that are connected to node $I$}
\label{fig:nodalpatch}
\end{figure}

Each node of the mesh is associated with their own patch of virtual elements.
The patch for node $I$ is denoted by $\mathcal{T}_I$ and is defined as the set of 
virtual elements connected to node $I$ (see~\fref{fig:nodalpatch}). Each node of a virtual element $E$ 
in the patch is assigned the area $\frac{1}{N_E^V}|E|$; that is,
the area of an element is uniformly distributed among its nodes. The representative 
area of node $I$ is denoted by $|I|$ and is computed by addition of all the areas 
that are assigned to node $I$ from the elements in $\mathcal{T}_I$; that is,
\begin{equation}\label{eq:nodal_area}
|I|=\sum_{E\in\mathcal{T}_I}\frac{1}{N_E^V}|E|.
\end{equation}
Similarly, each node of a virtual element $E$ is uniformly assigned the 
strain $\frac{1}{N_E^V}\widehat{\bsym{\varepsilon}}(\vm{u}_h)$. On considering each 
strain assigned to node $I$ from the elements in $\mathcal{T}_I$, 
the \textit{node-based uniform strain} is defined as follows:
\begin{equation}\label{eq:nodal_strain}
\widehat{\bsym{\varepsilon}}_I(\vm{u}_h)=
\frac{1}{|I|}\sum_{E\in\mathcal{T}_I}|E|\frac{1}{N_E^V}\widehat{\bsym{\varepsilon}}(\vm{u}_h).
\end{equation}
Since $\widehat{\bsym{\varepsilon}}(\vm{u}_h)$ is by definition given at the element level, then
from~\eref{eq:nodal_strain} the following nodal averaging operator is proposed:
\begin{equation}\label{eq:nodal_avg_operator}
\pi_I[\,\cdot\,]=\frac{1}{|I|}\sum_{E\in\mathcal{T}_I}|E|\frac{1}{N_E^V}[\,\cdot\,]_{{}_E},
\end{equation}
where $[\,\cdot\,]_{{}_E}$ denotes evaluation over the element $E$.

\subsection{NVEM nodal stiffness matrix and nodal force vector}
\label{sec:nodstiffness}

The NVEM nodal stiffness matrix is developed by substituting the discretizations~\eref{eq:discrete_field},
the projection operator~\eref{eq:operator} and the nodal averaging operator~\eref{eq:nodal_avg_operator}
into the left-hand side of \eref{eq:nodal_vemweakform_disc} for a node $I$. This gives
\begin{equation}\label{eq:nvem_nodal_stiff}
a_I \bigl(\pi_I[\Pi\vm{u}_h],\pi_I[\Pi\vm{v}_h]\bigr) 
+ s_I\bigl(\pi_I[\vm{u}_h-\Pi\vm{u}_h],\pi_I[\vm{v}_h-\Pi\vm{v}_h]\bigr) 
= \vm{q}^\transpose \vm{K}_I \vm{d},
\end{equation}
where $\vm{d}$ and $\vm{q}$ are column vectors of element nodal displacements 
and element nodal values associated with $\vm{v}_h$, respectively; $\vm{K}_I$
is the NVEM nodal stiffness matrix given by
\begin{equation}\label{eq:nvem_nodal_stiff_matrix}
\vm{K}_I = \vm{K}_I^\cons + \vm{K}_I^\stab,\quad
\vm{K}_I^\cons=|I|\,\vm{B}_I^\transpose\,\vm{D}\,\vm{B}_I,\quad
\vm{K}_I^\stab = (\vm{I}-\vm{P})_I^\transpose\,\vm{S}\,(\vm{I}-\vm{P})_I,
\end{equation}
where $\vm{D}$ is the constitutive matrix and $\vm{S}$ is the stability matrix
both defined in Section~\ref{sec:stabilization};
$\vm{B}_I=\pi_I[\,\vm{B}\,]$ and $(\vm{I}-\vm{P})_I=\pi_I[\,\vm{I}-\vm{P}\,]$,
where $\vm{I}$ is the identity $(2N_E^V \times 2N_E^V)$ matrix, 
and $\vm{B}$ and $\vm{P}$ are defined as
\begin{equation}\label{eq:matrix_B}
\vm{B}=\smat{
\vm{B}_1 & \cdots & \vm{B}_a & \cdots & \vm{B}_{N_E^V}}, \quad
\vm{B}_a = \smat{q_{1a} & 0 \\ 0 & q_{2a} \\ q_{2a} & q_{1a}},
\end{equation}
where $q_{ia}=\frac{1}{|E|}\int_{\partial E}\phi_a(\vm{x})n_i\,\diffs$
and can be exactly computed on the element boundary using a trapezoidal
rule giving the following algebraic expression:
\begin{equation}\label{eq:qia}
q_{ia}=\frac{1}{2|E|}\left(|e_{a-1}| n_{i(a-1)}+|e_a| n_{ia}\right),\quad i=1,2,
\end{equation}
where $n_{ia}$ is the $i$-th component of $\vm{n}_a$ and $|e_a|$ is the length of
the edge incident to node $a$ as defined in \fref{fig:polyelement};
\begin{equation}\label{eq:pimatrix}
\vm{P}=\vm{H}\vm{B}+\vm{G}\vm{R},
\end{equation}
where
\begin{equation}\label{eq:matrix_H}
\vm{H}=\smat{\vm{H}_1 & \cdots & \vm{H}_a & \cdots & \vm{H}_{N_E^V}}^\transpose, \quad
\vm{H}_a = \smat{(x_{1a}-\bar{x}_1) & 0 \\ 0 & (x_{2a}-\bar{x}_2) \\ \frac{1}{2}(x_{2a}-\bar{x}_2) & \frac{1}{2}(x_{1a}-\bar{x}_1)}^\transpose;
\end{equation}
\begin{equation}\label{eq:matrix_G}
\vm{G}=\smat{
\vm{G}_1 & \cdots & \vm{G}_a & \cdots & \vm{G}_{N_E^V}}^\transpose, \quad
\vm{G}_a = \smat{1 & 0 \\ 0 & 1 \\ \frac{1}{2}(x_{2a}-\bar{x}_2) & \frac{1}{2}(\bar{x}_1-x_{1a})}^\transpose;
\end{equation}
and
\begin{equation}\label{eq:matrix_R}
\vm{R}=\smat{
\vm{R}_1 & \cdots & \vm{R}_a & \cdots & \vm{R}_{N_E^V}}, \quad
\vm{R}_a = \smat{\frac{1}{N_E^V} & 0 \\ 0 & \frac{1}{N_E^V} \\ q_{2a} & -q_{1a}}.
\end{equation}

Similarly, the NVEM nodal force vector is developed using the right-hand side 
of~\eref{eq:nodal_vemweakform_disc} for a node $I$, which leads to
\begin{equation}\label{eq:nvem_nodforce}
 \ell_I\bigl(\pi_I[\Pi\vm{v}_h]\bigr) = \vm{q}^\transpose \vm{f}_I,
\end{equation}
where $\vm{q}$ is the column vector of element nodal values associated with $\vm{v}_h$,
and $\vm{f}_I$ is the NVEM nodal force vector associated with the 
body force and external tractions defined as
\begin{equation}\label{eq:nvem_nodal_forcevector}
\vm{f}_I = \vm{f}_I^b + \vm{f}_I^t,\quad
\vm{f}_I^b=|I|\,\bar{\vm{N}}_I^\transpose\widehat{\vm{b}}_I,\quad
\vm{f}_I^t = |I_\Gamma|\,\bar{\vm{N}}_{\Gamma,I}^\transpose\widehat{\vm{t}}_{N,I}.
\end{equation}

For computing the nodal body force vector $\vm{f}_I^b$ in~\eref{eq:nvem_nodal_forcevector}, 
$\bar{\vm{N}}_I=\pi_I[\,\bar{\vm{N}}\,]$ and $\widehat{\vm{b}}_I=\pi_I[\,\widehat{\vm{b}}\,]$,
where
\begin{equation}\label{eq:bar_N}
\bar{\vm{N}} = \smat{\bar{\vm{N}}_1 & \cdots & \bar{\vm{N}}_a & \cdots & \bar{\vm{N}}_{N_E^V}},\quad
\bar{\vm{N}}_a = \smat{\frac{1}{N_E^V} & 0 \\ 0 & \frac{1}{N_E^V}};
\end{equation}
and
\begin{equation}\label{eq:hat_bodyforce}
\widehat{\vm{b}}=\frac{1}{|E|}\int_E\vm{b}\,\diffx.
\end{equation}

Regarding the nodal traction force vector $\vm{f}_I^t$ in~\eref{eq:nvem_nodal_forcevector},
the nodal components are now computed with respect to the one-dimensional domain
on the Neumann boundary; that is, the representative nodal area reduces to a representative
nodal length $|I_\Gamma|=\sum_{e\in\mathcal{T}_I}\frac{1}{2}|e|$, where $e$ is
an element's edge located on the Neumann boundary and $|e|$ its length; $\mathcal{T}_I$ now 
represents the set of edges connected to node $I$ on the Neumann boundary.
Using the preceding definitions, the nodal averaging operator on the Neumann boundary 
is defined as
\begin{equation}\label{eq:nodal_avg_operator_neumann}
\pi_{I,\Gamma}[\,\cdot\,]=\frac{1}{|I_\Gamma|}\sum_{e\in\mathcal{T}_I}|e|\frac{1}{2}[\,\cdot\,]_e,
\end{equation}
where $[\,\cdot\,]_e$ denotes evaluation over the edge $e$. The remainder nodal matrices 
are then obtained as $\bar{\vm{N}}_{\Gamma,I}=\pi_{I,\Gamma}\bigl[\bar{\vm{N}}_\Gamma\bigr]$
and $\widehat{\vm{t}}_{N,I}=\pi_{I,\Gamma}\bigl[\,\widehat{\vm{t}}_{N}\bigr]$, where
\begin{equation}\label{eq:nodal_N_gamma}
\bar{\vm{N}}_\Gamma =
\smat{\frac{1}{2} & 0 & \frac{1}{2} & 0 \\ 0 & \frac{1}{2} & 0 & \frac{1}{2}},
\end{equation}
and
\begin{equation}\label{eq:nodal_traction}
\widehat{\vm{t}}_{N}=\frac{1}{|e|}\int_e\vm{t}_N\,\diffs.
\end{equation}

\subsection{NVEM equilibrium equations for elastoplasticity}
\label{sec:equilibrium}

For linear elastostatics as developed in Ref.~\cite{Ortiz-Silva-Salinas-Hitschfeld-Luza-Rebolledo:2023}, the discrete equilibrium
equation is obtained by the discretization of \eref{eq:nodal_vemweakform_disc}.
This is accomplished by summing \eref{eq:nvem_nodal_stiff} and \eref{eq:nvem_nodforce}
through all the nodes in the domain and invoking the arbitrariness of $\vm{q}$.
This results in the following system of equations:
\begin{equation}\label{eq:elastic_equilibrium}
\sum_{I\in\mathcal{T}_h} \biggl[|I|\,\vm{B}_I^\transpose\,\vm{D}\,\vm{B}_I
+ (\vm{I}-\vm{P})_I^\transpose\,\vm{S}\,(\vm{I}-\vm{P})_I\biggr]\vm{d} = \sum_{I\in\mathcal{T}_h} \vm{f}_I.
\end{equation}

Eq.~\eref{eq:elastic_equilibrium} is also the discrete equilibrium equation 
for elastoplasticity when $\vm{D}$ and $\vm{S}$ are nonlinear functions obtained
from the elastoplastic constitutive law. As usual, the solution for this problem
requires the linearization of~\eref{eq:elastic_equilibrium}. Doing this gives
\begin{equation}\label{eq:elastoplastic_equilibrium}
\begin{split}
\sum_{I\in\mathcal{T}_h}&\biggl[|I|\,\vm{B}_I^\transpose\,\breve{\vm{D}}^{ep}\,\vm{B}_I
+ (\vm{I}-\vm{P})_I^\transpose\,\breve{\vm{S}}\,(\vm{I}-\vm{P})_I\biggr]\Delta\vm{d} \\
&= - \sum_{I\in\mathcal{T}_h} \biggl[|I|\,\vm{B}_I^\transpose\,\breve{\bsym{\sigma}}
+ (\vm{I}-\vm{P})_I^\transpose\,\breve{\vm{S}}\,(\vm{I}-\vm{P})_I\,\vm{d}- \vm{f}_I\biggr],
\end{split}
\end{equation}
where $\breve{\vm{D}}^{ep}$, $\breve{\vm{S}}$, and $\breve{\bsym{\sigma}}$,
which are given in Section~\ref{sec:stabilization}, are the elastoplastic 
consistent tangent operator, the stability matrix, and the nonlinear stress, respectively, 
evaluated at node $I$ (by means of the elastoplastic constitutive law) using 
the node-based uniform strain $\widehat{\bsym{\varepsilon}}_I$.

The linearized equilibrium equation~\eref{eq:elastoplastic_equilibrium} is
used to solve the equilibrium state 
$\vm{d}_{n+1}^{(k)}=\vm{d}_{n+1}^{(k-1)}+\Delta\vm{d}^{(k)}$ 
at time $t_{n+1}$ with a time increment $\Delta t = t_{n+1}-t_n$
via Newton-Raphson iterations, as follows:
\begin{equation}\label{eq:nr_scheme}
\sum_{I\in\mathcal{T}_h}  \biggl[\vm{K}_{I,T}^\cons+\vm{K}_{I,T}^\stab\biggr]_{n+1}^{(k-1)}\Delta\vm{d}^{(k)}
= -\sum_{I\in\mathcal{T}_h,}\biggl[\vm{f}_I^\cons \, + \, \vm{f}_I^\stab -\vm{f}_I\biggr]_{n+1}^{(k-1)},
\end{equation}
where
\begin{equation}\label{eq:nr_scheme2}
\begin{split}
\vm{K}_{I,T}^\cons &=|I|\,\vm{B}_I^\transpose\,\breve{\vm{D}}^{ep}\,\vm{B}_I,\quad
\vm{K}_{I,T}^\stab =(\vm{I}-\vm{P})_I^\transpose\,\breve{\vm{S}}\,(\vm{I}-\vm{P})_I,\\ 
&\vm{f}_I^\cons=|I|\,\vm{B}_I^\transpose\,\breve{\bsym{\sigma}},\quad
\vm{f}_I^\stab=(\vm{I}-\vm{P})_I^\transpose\,\breve{\vm{S}}\,(\vm{I}-\vm{P})_I\,\vm{d}.
\end{split}
\end{equation}

\section{Elastoplastic constitutive model and stabilization}
\label{sec:stabilization}

Within the standard VEM framework, stabilization is one of the key ingredients for 
convergence of the method. However, in nodal integration, stabilization can make the formulation 
somewhat stiff in incompressible settings~\cite{Puso-Solberg:2006}. 
Therefore, the stabilization in any nodal integration scheme, \alejandro{which includes the NVEM, 
must be cautiously chosen to not jeopardize its locking-free essence. 
To deal with this issue, we propose a \textit{D-recipe} stabilization~\cite{BeiraoDaVeiga-Dassi-Russo:2017,Mascotto:2018}
that uses only a deviatoric term.} Firstly, a summary of the constitutive model considered 
for the elastoplastic solid is given and, secondly, the stabilization for the NVEM is detailed.

\subsection{Constitutive model}
\label{sec:constitutivemodel}

The strain tensor $\bsym{\varepsilon}$ is split into elastic ($\bsym{\varepsilon}^e$) 
and plastic ($\bsym{\varepsilon}^p$) parts; that is,
\begin{equation}\label{eq:strainsplit}
\bsym{\varepsilon} = \bsym{\varepsilon}^e+\bsym{\varepsilon}^p.
\end{equation}
The elastic part is governed by the standard linear elastic law, and
the plastic part by the von Mises model with mixed linear hardening 
(for details on this model, see for instance Ref.~\cite{deSouza-Peric-Owen:2008}). 
The yield function for this model is
\begin{equation}\label{eq:claw}
\Phi(\bsym{\sigma},\bsym{\beta},\sigma_y) = 
\sqrt{3\,J_2(\vm{s}(\bsym{\sigma})-\bsym{\beta})}-\sigma_y(\bar{\varepsilon}^p)
= \sqrt{\frac{3}{2}} \|\bsym{\eta}\|-\sigma_y(\bar{\varepsilon}^p),
\end{equation}
where $\bsym{\beta}$ is the backstress tensor, $\sigma_y$ is a function
of the accumulated plastic strain $\bar{\varepsilon}^p$
and defines the radius of the yield surface, and $\bsym{\eta}$ is the
relative stress given by
\begin{equation}\label{eq:relstress}
\bsym{\eta}=\vm{s}-\bsym{\beta},
\end{equation}
where $\vm{s}$ is the deviatoric stress. The plastic flow is described by
the following associative law:
\begin{equation}\label{eq:plasticflow}
\dot{\bsym{\varepsilon}}^p=\dot{\gamma}\sqrt{\frac{3}{2}}\frac{\bsym{\eta}}{\|\bsym{\eta}\|},
\end{equation}
where $\dot{\bsym{\varepsilon}}^p$ is the rate of the plastic strain tensor
and $\dot{\gamma}$ is the plastic multiplier. 

The mixed linear hardening combines
linear isotropic and linear kinematic hardening models. The linear isotropic hardening
model is defined by the linear function
\begin{equation}\label{eq:linisohardening}
\sigma_y(\bar{\varepsilon}^p)=\sigma_{y0}+H_i\,\bar{\varepsilon}^p,
\end{equation}
where $\sigma_{y0}$ is the initial yield stress and $H_i$ is the linear
isotropic hardening modulus. The linear kinematic hardening model
describes the evolution law for the backstress as the following
linear function:
\begin{equation}\label{eq:linkinhardening}
\dot{\bsym{\beta}}=\frac{2}{3}H_k\,\dot{\bsym{\varepsilon}}^p,
\end{equation}
where $H_k$ is the linear kinematic hardening modulus.

Let $n$ and $n+1$ be the subindices that denote the previous and current states, respectively.
These subindices are used for labeling the time at which the variables
that are involved in the constitutive law are evaluated. Using the preceding notation, 
the elastoplastic consistent tangent operator (in Voigt notation) for the above model is
\begin{equation}\label{eq:eptangent}
\begin{split}
\vm{D}_{n+1}^{ep} ={} & 2G\,\Biggl(1-\frac{\Delta\gamma 3G}{\bar{q}_{n+1}^\mathrm{trial}}\Biggr)\,\vm{I_d}
+ 6G^2\Biggl(\frac{\Delta\gamma}{\bar{q}_{n+1}^\mathrm{trial}}
-\frac{1}{3G+H_k+H_i}\Biggr)\,\bar{\vm{N}}_{n+1}\bar{\vm{N}}_{n+1}^\transpose\\
& + K\vm{m}\vm{m}^\transpose,
\end{split}
\end{equation}
where $G$ and $K$ are the shear modulus and bulk modulus of the material, respectively.
Here, we use the three-dimensional constitutive law from where the plain strain state components
are extracted. The remainder quantities that appear in~\eref{eq:eptangent} are defined as follows:
\begin{equation}\label{eq:idmatrix}
\vm{I}_d = \smat{1&0&0&0&0&0\\0&1&0&0&0&0\\0&0&1&0&0&0\\0&0&0&0.5&0&0\\0&0&0&0&0.5&0\\0&0&0&0&0&0.5}-\frac{1}{3}\vm{m}\vm{m}^\transpose,\quad
\vm{m} = \smat{1&1&1&0&0&0}^\transpose;
\end{equation}
\begin{equation}\label{eq:qntrial}
\begin{split}
&\bar{q}_{n+1}^\mathrm{trial}=\sqrt{\frac{3}{2}}\|\bsym{\eta}_{n+1}^\mathrm{trial}\|,\quad
\bsym{\eta}_{n+1}^\mathrm{trial}=\vm{s}_{n+1}^\mathrm{trial}-\bsym{\beta}_n,\\
\vm{s}_{n+1}^\mathrm{trial} =2G &\,\Bigl(\bsym{\varepsilon}_{n+1}^{e\,\,\mathrm{trial}} 
- \frac{1}{3}\trace{\bsym{\varepsilon}_{n+1}^{e\,\,\mathrm{trial}}}\,\vm{m}\Bigr),\quad
\bsym{\varepsilon}_{n+1}^{e\,\,\mathrm{trial}}=\bsym{\varepsilon}_{n+1}-\bsym{\varepsilon}_n^p;
\end{split}
\end{equation}
\begin{equation}\label{eq:deltagamma}
\Delta\gamma = \frac{\Phi^\mathrm{trial}}{3G+H_k+H_i},\quad
\Phi^\mathrm{trial}=\bar{q}_{n+1}^\mathrm{trial}-
\sigma_y(\bar{\varepsilon}_{n+1}^{p\,\,\mathrm{trial}}),\quad
\bar{\varepsilon}_{n+1}^{p\,\,\mathrm{trial}}=\bar{\varepsilon}_n^p;
\end{equation}
\begin{equation}\label{eq:nbarvector}
\bar{\vm{N}}_{n+1}=\frac{\bsym{\eta}_{n+1}^\mathrm{trial}}{\|\bsym{\eta}_{n+1}^\mathrm{trial}\|}.
\end{equation}

The strain/stress update is performed as follows:
\begin{equation}\label{eq:stressupdate}
\begin{split}
\bar{\varepsilon}_{n+1}^p &= \bar{\varepsilon}_n^p+\Delta\gamma,\quad
\bsym{\varepsilon}_{n+1}^p = \bsym{\varepsilon}_n^p+\Delta\gamma\sqrt{\frac{3}{2}}\,\bar{\vm{N}}_{n+1},\\
\bsym{\beta}_{n+1} &= \bsym{\beta}_n + \Delta\gamma\sqrt{\frac{2}{3}}H_k\bar{\vm{N}}_{n+1},\\
\vm{s}_{n+1} &=\vm{s}_{n+1}^{\mathrm{trial}}-2G\Delta\gamma\sqrt{\frac{3}{2}}\,\bar{\vm{N}}_{n+1},\\
\bsym{\sigma}_{n+1} &=\vm{s}_{n+1}+K\,\trace{\bsym{\varepsilon}_{n+1}^{e\,\,\mathrm{trial}}}\,\vm{m}.
\end{split}
\end{equation}

As previously mentioned in Section~\ref{sec:equilibrium}, 
in the NVEM the constitutive law is
evaluated using the nodal strain $\widehat{\bsym{\varepsilon}}_I$;
that is, in~\eref{eq:nr_scheme2},
$\breve{\vm{D}}^{ep}=\vm{D}^{ep}(\bsym{\varepsilon}=\widehat{\bsym{\varepsilon}}_I)$
and 
$\breve{\bsym{\sigma}}=\bsym{\sigma}(\bsym{\varepsilon}=\widehat{\bsym{\varepsilon}}_I)$.

The above computations of the elastoplastic consistent tangent operator
and strain/stress update are part of the implicit elastic predictor/return mapping 
algorithm for numerical integration of the von Mises model with mixed linear hardening~\cite{deSouza-Peric-Owen:2008}. 
This algorithm is summarized in Algorithm~\ref{algo1}.

\begin{algorithm}
\setstretch{1.4}
\caption{Implicit elastic predictor/return mapping algorithm for the von Mises model with mixed linear hardening}\label{algo1}
\begin{algorithmic}[]
\Require $\bigl(\widehat{\bsym{\varepsilon}}_I\bigr)_{n+1}$, $\bar{\varepsilon}_n^p$,
        $\bsym{\varepsilon}_n^p$, $\bsym{\beta}_n$
\Ensure $\bar{\varepsilon}_{n+1}^p$, $\bsym{\varepsilon}_{n+1}^p$, $\bsym{\beta}_{n+1}$, $\breve{\bsym{\sigma}}_{n+1}$, $\breve{\vm{D}}_{n+1}^{ep}$ 
\medskip
\State (Elastic predictor)
\State Set $\bsym{\varepsilon}_{n+1} = \bigl(\widehat{\bsym{\varepsilon}}_I\bigr)_{n+1}$
\State $\bsym{\varepsilon}_{n+1}^{e\,\,\mathrm{trial}}=\bsym{\varepsilon}_{n+1}-\bsym{\varepsilon}_n^p$
\State $\vm{s}_{n+1}^\mathrm{trial}=2G\,\Bigl(\bsym{\varepsilon}_{n+1}^{e\,\,\mathrm{trial}} 
            - \frac{1}{3}\trace{\bsym{\varepsilon}_{n+1}^{e\,\,\mathrm{trial}}}\,\vm{m}\Bigr)$
\State $\bsym{\eta}_{n+1}^\mathrm{trial}=\vm{s}_{n+1}^\mathrm{trial}-\bsym{\beta}_n, \quad
            \bar{q}_{n+1}^\mathrm{trial}=\sqrt{\frac{3}{2}}\|\bsym{\eta}_{n+1}^\mathrm{trial}\|$
\State $\bar{\varepsilon}_{n+1}^{p\,\,\mathrm{trial}}=\bar{\varepsilon}_n^p, \quad
            \bsym{\beta}_{n+1}^\mathrm{trial}=\bsym{\beta}_n$
\medskip
\State (Elastic/plastic check and update state)
\State $\Phi^\mathrm{trial}=\bar{q}_{n+1}^\mathrm{trial}-\sigma_y(\bar{\varepsilon}_{n+1}^{p\,\,\mathrm{trial}})$
\If{$\Phi^\mathrm{trial}\leq 0$}\label{algln2} \quad(elastic state)
            \State $\bar{\varepsilon}_{n+1}^p = \bar{\varepsilon}_{n+1}^{p\,\,\mathrm{trial}},
            \quad \bsym{\varepsilon}_{n+1}^p=\bsym{\varepsilon}_n^p$
            \State $\bsym{\beta}_{n+1}=\bsym{\beta}_{n+1}^\mathrm{trial}$
            \State $\bsym{\sigma}_{n+1}=\vm{s}_{n+1}^\mathrm{trial}+K\,\trace{\bsym{\varepsilon}_{n+1}^{e\,\,\mathrm{trial}}}\,\vm{m}$
            \State $\vm{D}_{n+1}^{ep} = 2G\vm{I}_d + K\vm{m}\vm{m}^\transpose$
\Else \quad(plastic corrector)
            \State $\Delta\gamma = \frac{\Phi^\mathrm{trial}}{3G+H_k+H_i}, \quad 
            \bar{\vm{N}}_{n+1}=\frac{\bsym{\eta}_{n+1}^\mathrm{trial}}{\|\bsym{\eta}_{n+1}^\mathrm{trial}\|}$  
            \State $\bar{\varepsilon}_{n+1}^p = \bar{\varepsilon}_n^p+\Delta\gamma, \quad                  
            \bsym{\varepsilon}_{n+1}^p = \bsym{\varepsilon}_n^p+\Delta\gamma\sqrt{\frac{3}{2}}\,\bar{\vm{N}}_{n+1}$
            \State $\bsym{\beta}_{n+1} = \bsym{\beta}_n + \Delta\gamma\sqrt{\frac{2}{3}}H_k\,\bar{\vm{N}}_{n+1}$
            \State $\vm{s}_{n+1}=\vm{s}_{n+1}^{\mathrm{trial}}-2G\Delta\gamma\sqrt{\frac{3}{2}}\,\bar{\vm{N}}_{n+1}$
            \State $\bsym{\sigma}_{n+1}=\vm{s}_{n+1}+K\,\trace{\bsym{\varepsilon}_{n+1}^{e\,\,\mathrm{trial}}}\,\vm{m}$
            \State $\vm{D}_{n+1}^{ep}=2G\,\Bigl(1-\frac{\Delta\gamma 3G}{\bar{q}_{n+1}^\mathrm{trial}}\Bigr)\,\vm{I_d} 
            + 6G^2\Bigl(\frac{\Delta\gamma}{\bar{q}_{n+1}^\mathrm{trial}} 
            -\frac{1}{3G+H_k+H_i}\Bigr)\,\bar{\vm{N}}_{n+1}\bar{\vm{N}}_{n+1}^\transpose 
            + K\vm{m}\vm{m}^\transpose$
\EndIf
\State $\breve{\bsym{\sigma}}_{n+1}=\bsym{\sigma}_{n+1}$, $\breve{\vm{D}}_{n+1}^{ep}=\vm{D}_{n+1}^{ep}$
\end{algorithmic}
\end{algorithm}

\subsection{Stabilization}
\label{sec:stabmatrix}

\alejandro{As mentioned before, some precautions must be taken to stabilize the NVEM
so that the locking-free behavior of the nodal integration scheme is preserved.
Thus, drawing inspiration from one of the possibilities 
already explored in Ref.~\cite{Ortiz-Silva-Salinas-Hitschfeld-Luza-Rebolledo:2023},
the stabilization issue in the NVEM for elastoplastic solids is dealt with
a diagonal stability matrix that uses only a deviatoric term, as follows:}
\begin{equation}\label{eq:nvemstab}
\bigl(\breve{\vm{S}}\bigr)_{i,i}=
\max{\left(1,\Bigl[|I|\,\vm{B}_I^\transpose\,\vm{D}_\mathrm{d}^{e}\,\vm{B}_I\Bigr]_{i,i}\right)},
\end{equation}
where $\vm{D}_\mathrm{d}^{e}$ is the deviatoric part of the elastic moduli and is given by
\begin{equation}\label{eq:elasticdevmoduli}
\vm{D}_\mathrm{d}^{e}=2G\vm{I}_d.
\end{equation}
This choice performs very well in a variety of \alejandro{two-dimensional} numerical tests 
and does not introduce any tuning parameter.

\section{Numerical examples}
\label{sec:numexp}

In this section, some benchmark tests are conducted to demonstrate the performance of the
NVEM in elastoplastic solids simulations. The method is compared with the well-known
locking-free 9-node B-bar quadrilateral finite element~\cite{Simo-Hughes:1998} (FEM Q9 B-bar)
as well as with the standard linearly precise virtual element (VEM)
and in some cases with the standard 4-node quadrilateral finite element (FEM Q4).
All the tests are conducted using the constitutive model described 
in Section~\ref{sec:constitutivemodel}.
Throughout this section, DOF stands for degree(s) of freedom,
$E_\mathrm{Y}$ is the Young's modulus, $\nu$ is the Poisson's ratio, $\sigma_{y0}$ is
the initial yield stress, and $H_i$ and $H_k$ are the linear isotropic hardening
modulus and the linear kinematic hardening modulus, respectively.

\subsection{Thick-walled cylinder}
\label{sec:cylinder}

\begin{figure}[!htb]
\centering
\subfigure[] {\label{fig:twcylmeshes_a}\includegraphics[width=0.31\linewidth]{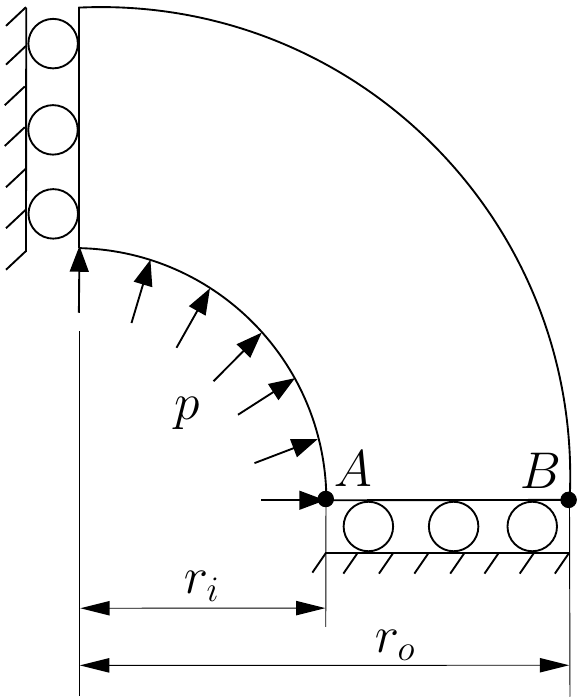}}
\subfigure[] {\label{fig:twcylmeshes_b}\includegraphics[width=0.4\linewidth]{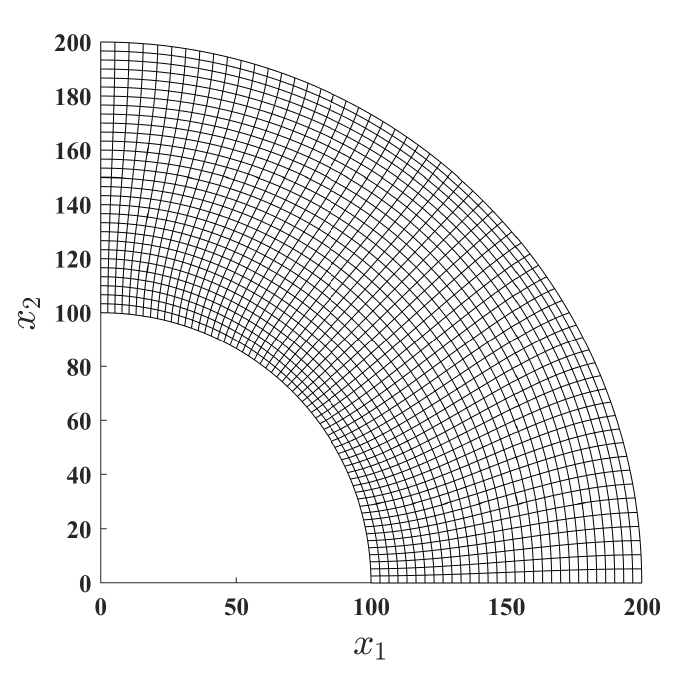}}
\caption{Thick-walled cylinder problem. (a) Geometry and boundary conditions, and 
        (b) mesh used for benchmarking the VEM, NVEM, FEM Q4 and FEM Q9 B-bar approaches}
\label{fig:twcyl_geomeshes}
\end{figure}

In this test, the ability of the NVEM for solving compressible and nearly incompressible 
elastoplastic problems is demonstrated. The problem consists of a plane strain 
(unit thickness) representation of a thick-walled cylinder under internal pressure. 
The geometry, boundary conditions, and mesh used in this numerical test are shown 
in~\fref{fig:twcyl_geomeshes}, where the internal pressure is $p=180$ MPa, and the 
internal and external radii are $r_i=100$ mm and $r_o=200$ mm, respectively. 
The material parameters used are the following: $E_\mathrm{Y}=210000$ MPa, $\nu=0.3$ 
for the compressible case and $\nu=0.4999$ for the nearly incompressible case,
$\sigma_{y0}=240$ MPa, $H_i = H_k =0$ MPa (perfect plasticity).

\begin{figure}[!htb]
\centering
\subfigure[]{\label{fig:twcyl_dispcurve_a}\includegraphics[width=0.49\linewidth]{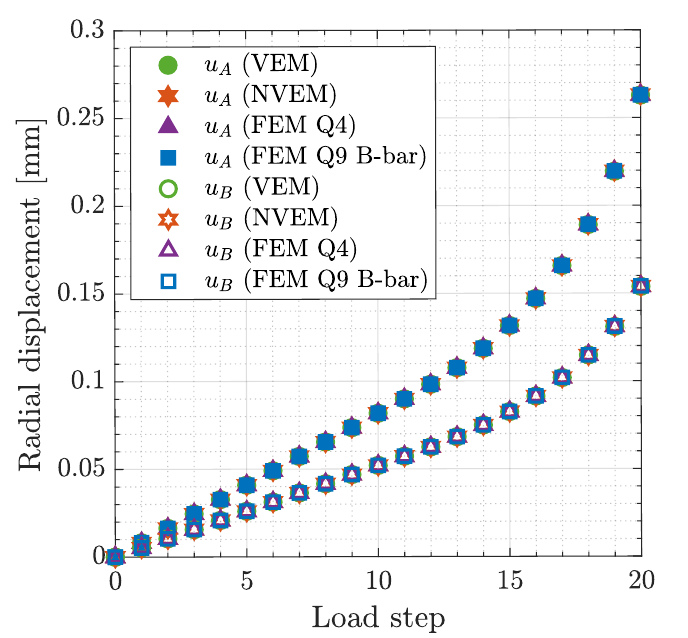}}
\subfigure[]{\label{fig:twcyl_dispcurve_b}\includegraphics[width=0.49\linewidth]{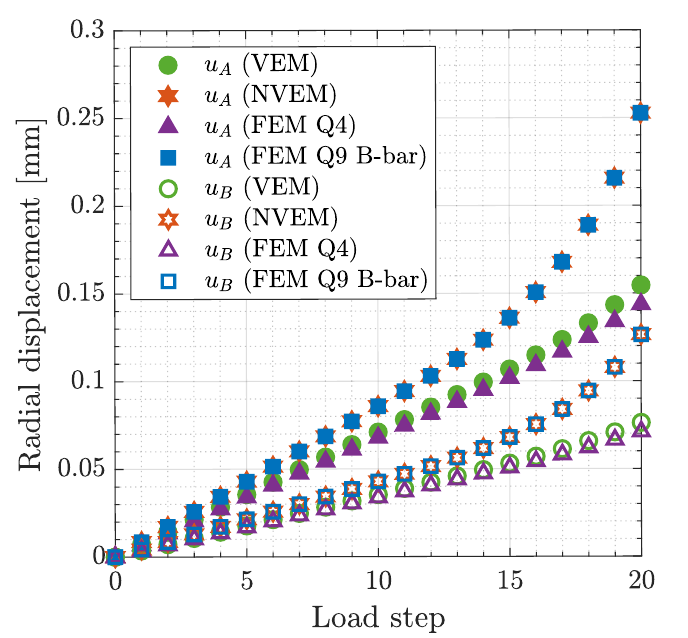}}
\caption{Radial displacement at points $A$ ($u_A$) and $B$ ($u_B$) 
         due to the applied internal pressure in steps for the thick-walled 
         cylinder problem. (a) Compressible case ($\nu=0.3$), and (b) nearly incompressible case ($\nu=0.4999$)}
\label{fig:twcyl_dispcurve}
\end{figure}

\begin{figure}[!tbp]
\centering
\mbox{
\subfigure[] {\label{fig:twcyl_disp_a}\includegraphics[width=0.49\linewidth]{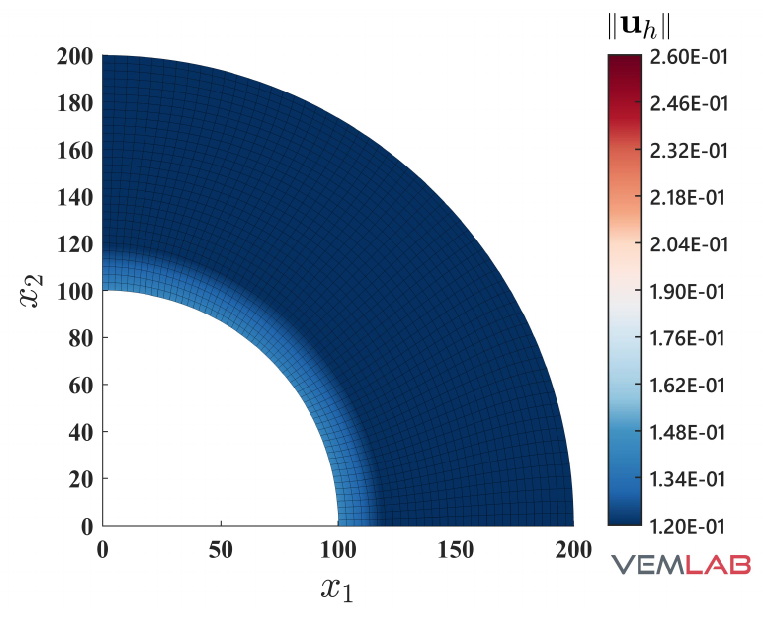}}
\subfigure[] {\label{fig:twcyl_disp_b}\includegraphics[width=0.49\linewidth]{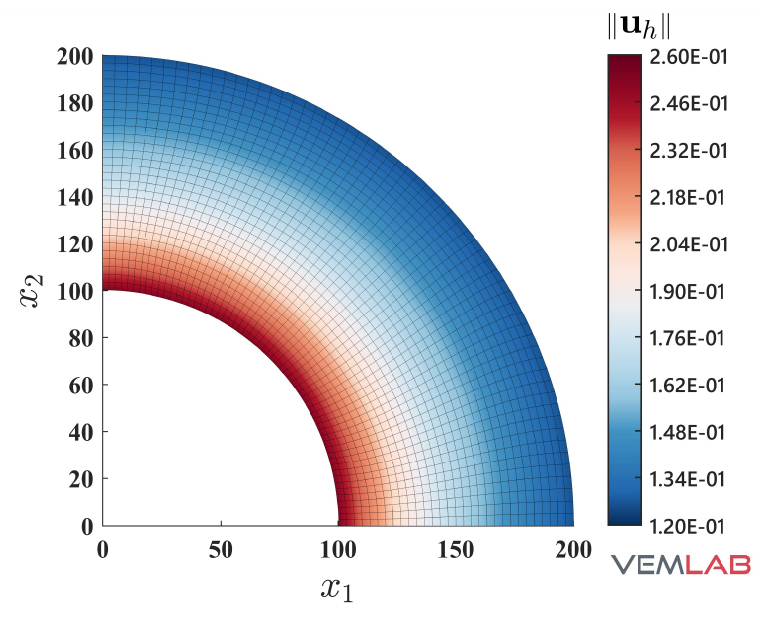}}
}
\mbox{
\subfigure[] {\label{fig:twcyl_disp_c}\includegraphics[width=0.49\linewidth]{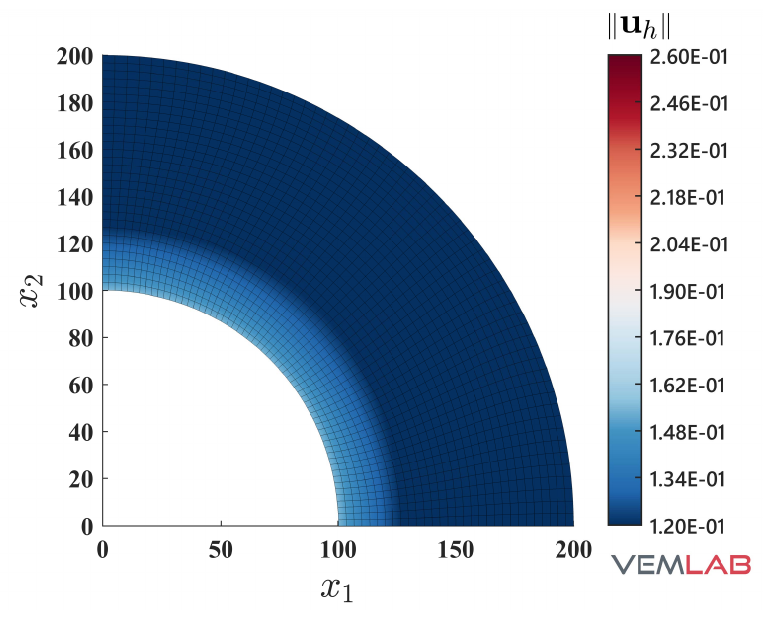}}
\subfigure[] {\label{fig:twcyl_disp_d}\includegraphics[width=0.49\linewidth]{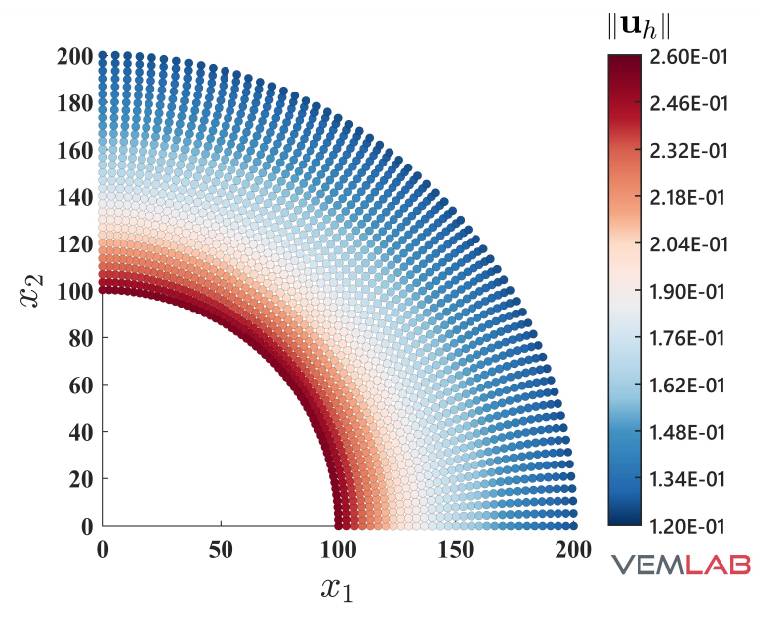}}
}
\caption{Plots of the total displacement field solution in mm at the last load step 
        for the nearly incompressible thick-walled cylinder problem ($\nu=0.4999$). 
        (a) FEM Q4, (b) FEM Q9 B-bar, (c) VEM, and (d) NVEM}
\label{fig:twcyl_totaldisp}
\end{figure}

\begin{figure}[!tbp]
\centering
\mbox{
\subfigure[] {\label{fig:twcyl_pressure_a}\includegraphics[width=0.49\linewidth]{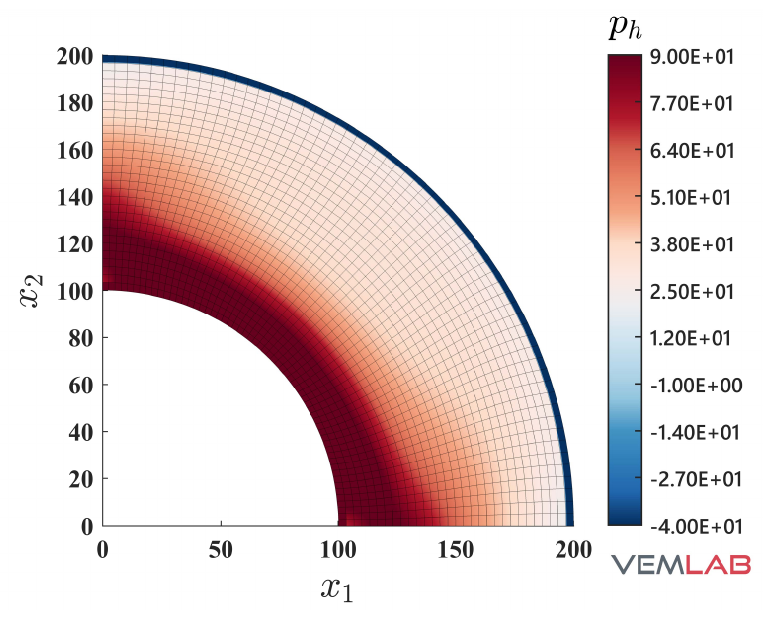}}
\subfigure[] {\label{fig:twcyl_pressure_b}\includegraphics[width=0.49\linewidth]{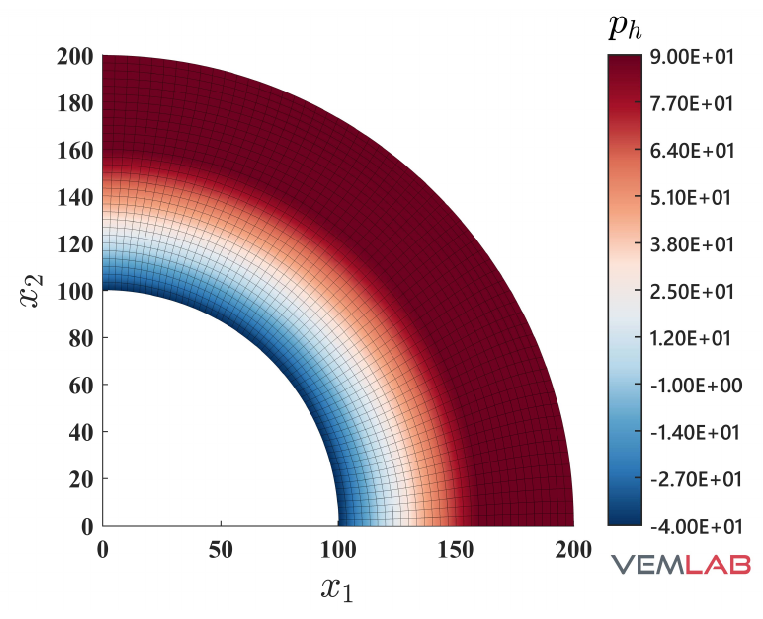}}
}
\mbox{
\subfigure[] {\label{fig:twcyl_pressure_c}\includegraphics[width=0.49\linewidth]{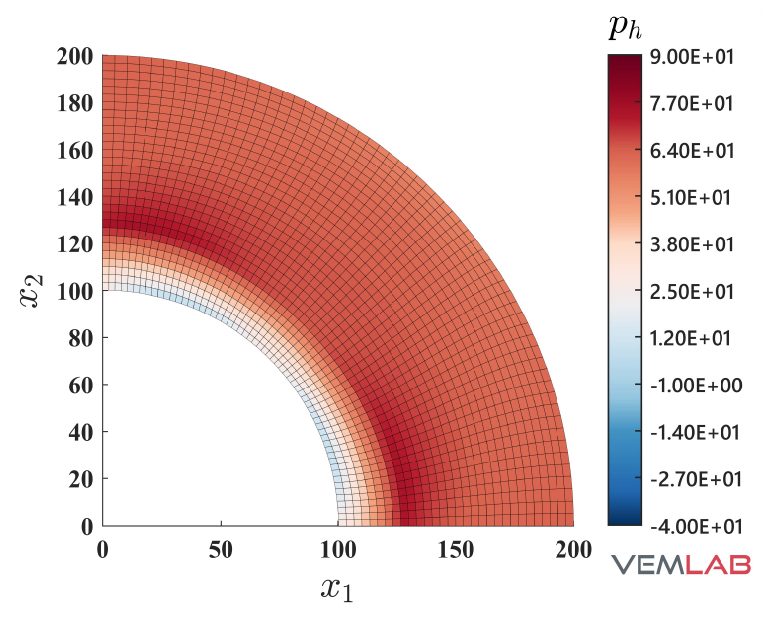}}
\subfigure[] {\label{fig:twcyl_pressure_d}\includegraphics[width=0.49\linewidth]{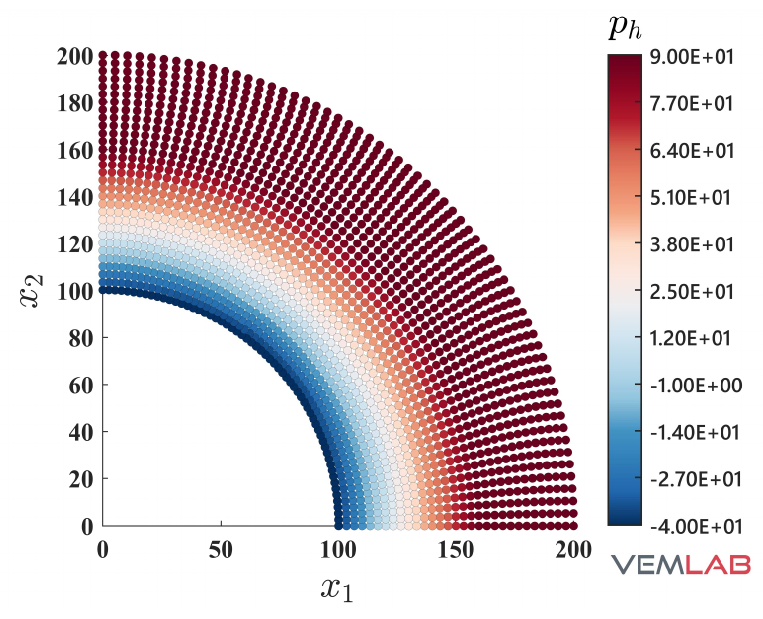}}
}
\caption{Plots of the pressure field solution in MPa for the nearly incompressible thick-walled cylinder problem ($\nu=0.4999$). 
        (a) FEM Q4, (b) FEM Q9 B-bar, (c) VEM, and (d) NVEM}
\label{fig:twcyl_pressure}
\end{figure}

\begin{figure}[!htb]
\centering
\mbox{
\subfigure[] {\label{fig:twcyl_vmstress_a}\includegraphics[width=0.49\linewidth]{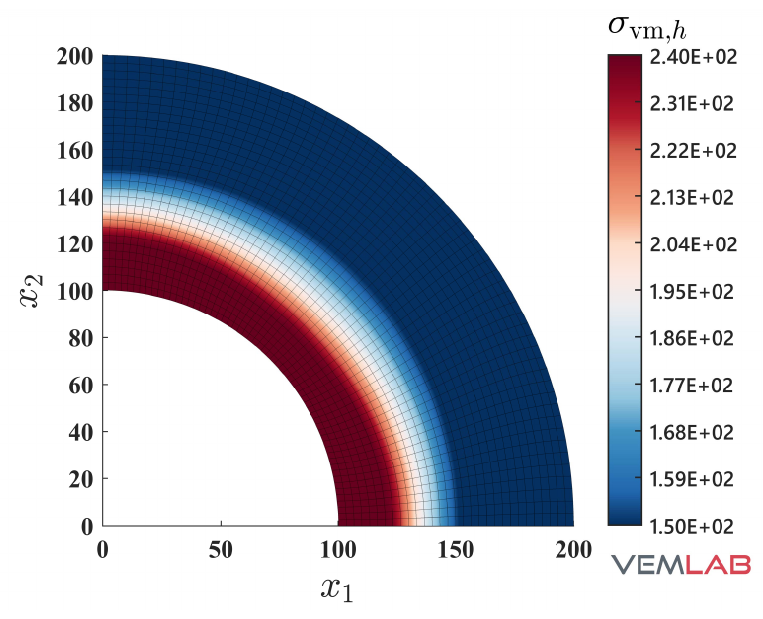}}
\subfigure[] {\label{fig:twcyl_vmstress_b}\includegraphics[width=0.49\linewidth]{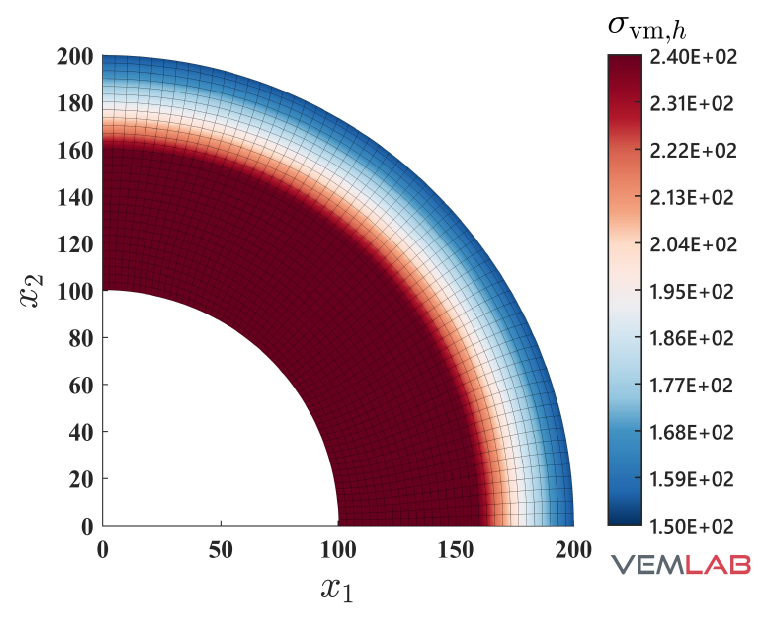}}
}
\mbox{
\subfigure[] {\label{fig:twcyl_vmstress_c}\includegraphics[width=0.49\linewidth]{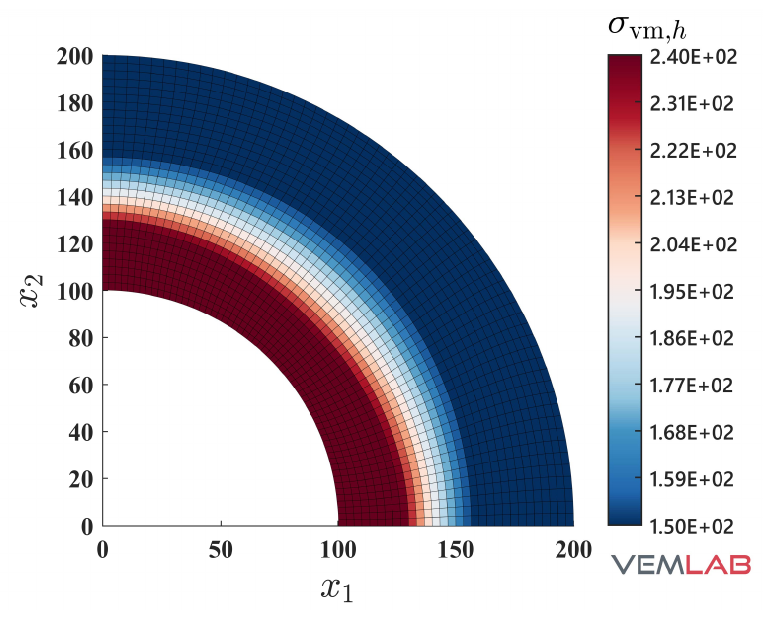}}
\subfigure[] {\label{fig:twcyl_vmstress_d}\includegraphics[width=0.49\linewidth]{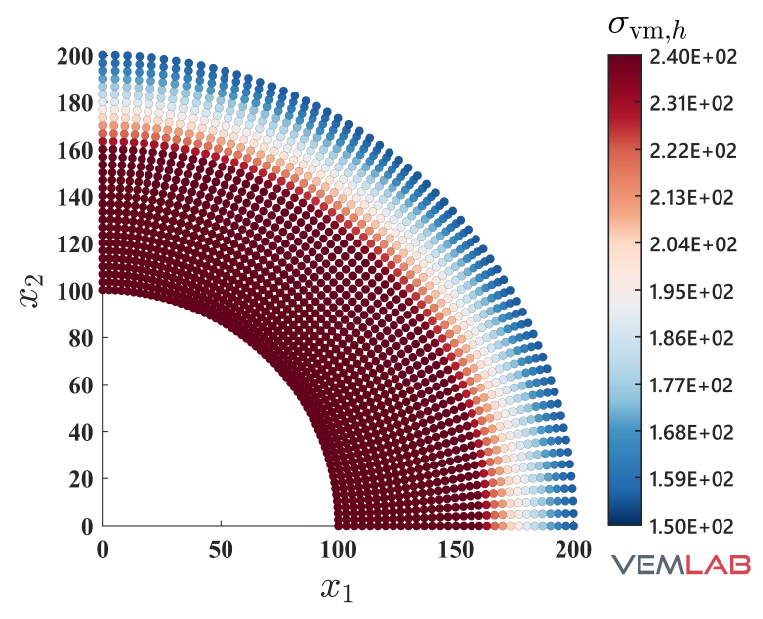}}
}
\caption{Plots of the von Mises stress field solution in MPa for the nearly incompressible thick-walled cylinder problem ($\nu=0.4999$). 
        (a) FEM Q4, (b) FEM Q9 B-bar, (c) VEM, and (d) NVEM}
\label{fig:twcyl_vmstress}
\end{figure}

The radial displacement at points $A$ and $B$ for the compressible and 
nearly incompressible cases is summarized in~\fref{fig:twcyl_dispcurve} for all the methods.
For the compressible case, all the methods match very well (\fref{fig:twcyl_dispcurve_a}). 
On the other hand, for the nearly incompressible case (\fref{fig:twcyl_dispcurve_b}), 
the NVEM and FEM Q9 B-bar methods still match very well, whereas the VEM and FEM Q4
clearly exhibit a locking behavior since radial displacements are smaller than expected.
The same behaviors are observed in the total displacement (\fref{fig:twcyl_totaldisp}), 
pressure (\fref{fig:twcyl_pressure}), and von Mises stress (\fref{fig:twcyl_vmstress}) 
field solutions at the last load step for the nearly incompressible case.

\subsection{Cook's membrane}
\label{sec:cook}

\begin{figure}[!htb]
\centering
\includegraphics[width=0.5\linewidth]{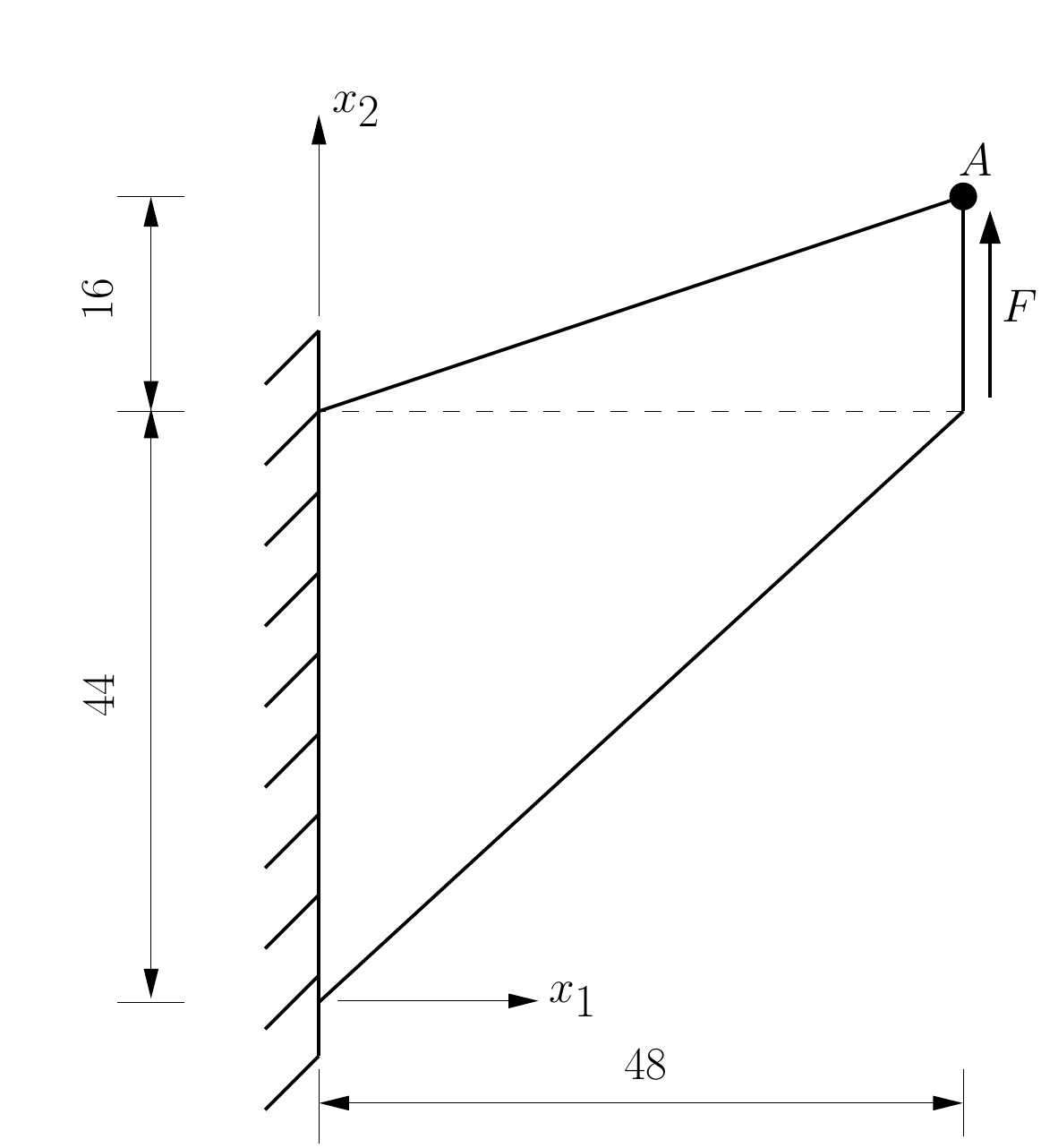}
\caption{Geometry and boundary conditions for the Cook's membrane problem with dimensions in mm}
\label{fig:cookproblem}
\end{figure}

\begin{figure}[!tbp]
\centering
\subfigure[] {\label{fig:cookmeshes_a}\includegraphics[width=0.49\linewidth]{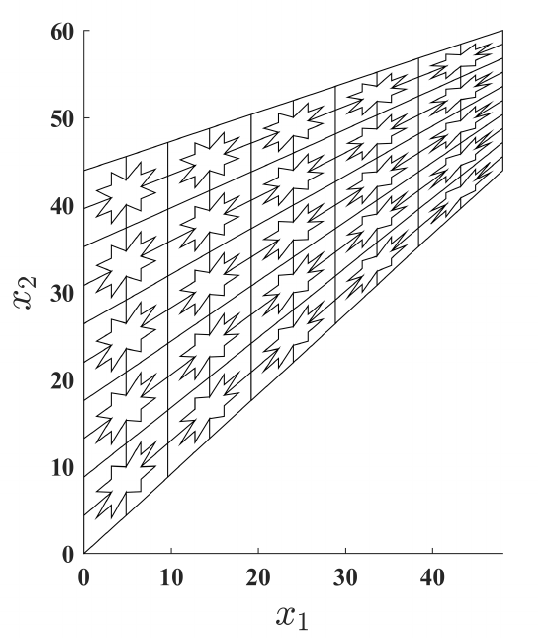}}
\subfigure[] {\label{fig:cookmeshes_b}\includegraphics[width=0.49\linewidth]{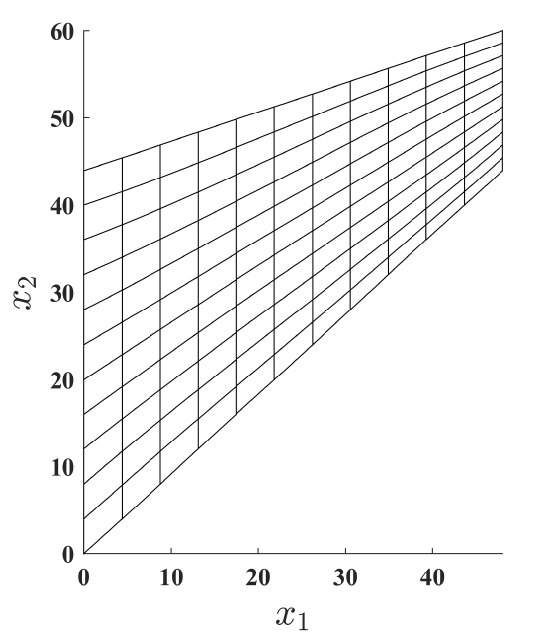}}
\caption{Sample meshes for the Cook's membrane problem that are used in the 
(a) VEM and NVEM, and (b) FEM Q9 B-bar approaches}
\label{fig:cookmeshes}
\end{figure}

The next example consists of a tapered beam fixed along one end and loaded with a shear force at
the other end. It is designed to study the performance of numerical formulations under combined
bending and shear when the solid material behaves nearly incompressible.
The geometry and boundary conditions are depicted in~\fref{fig:cookproblem}, where
the shear load is $F=3.6$ N/mm (total shear load of $57.6$ N). The beam has a unit thickness and plane strain
condition is assumed. The following material parameters are used:
$E_\mathrm{Y}=1500$ MPa, $\nu=0.4999$, $\sigma_{y0}=7.5$ MPa, $H_i=3.25$ MPa, and $H_k=0$ MPa.
In this test, the performance of the NVEM is compared with the 
VEM and the FEM Q9 B-bar. Sample meshes used in this test are shown
in \fref{fig:cookmeshes}.

\begin{figure}[!htb]
\centering
\includegraphics[width=0.6\linewidth]{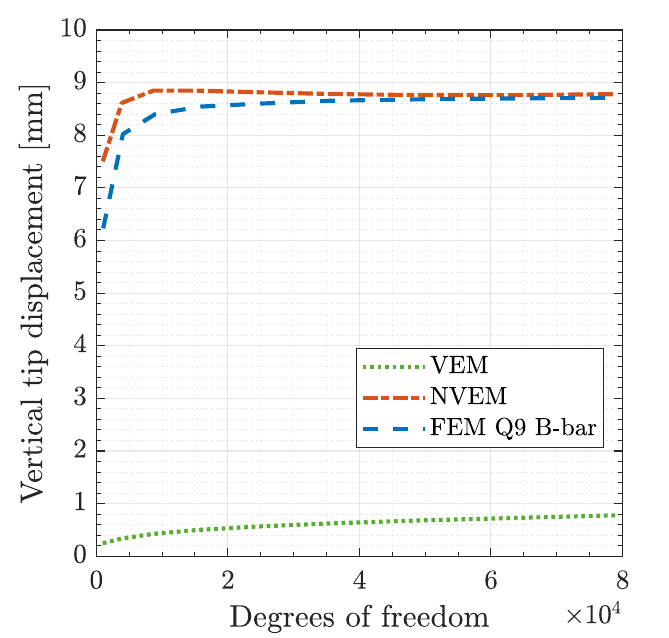}
\caption{Convergence of the vertical displacement at the tip of the Cook's membrane (point A) upon mesh refinement}
\label{fig:cookconvergence}
\end{figure}

The convergence of the vertical displacement at the tip of the beam
(point $A$ in~\fref{fig:cookproblem}) upon mesh refinement is depicted
in \fref{fig:cookconvergence}. As expected, the VEM solution exhibits
a severe locking behavior, whereas the NVEM and FEM Q9 B-bar are locking-free
and their solutions are in good agreement upon mesh refinement.

\begin{figure}[!tbp]
\centering
\mbox{
\subfigure[] {\label{fig:cook_solutions_a}\includegraphics[width=0.49\linewidth]{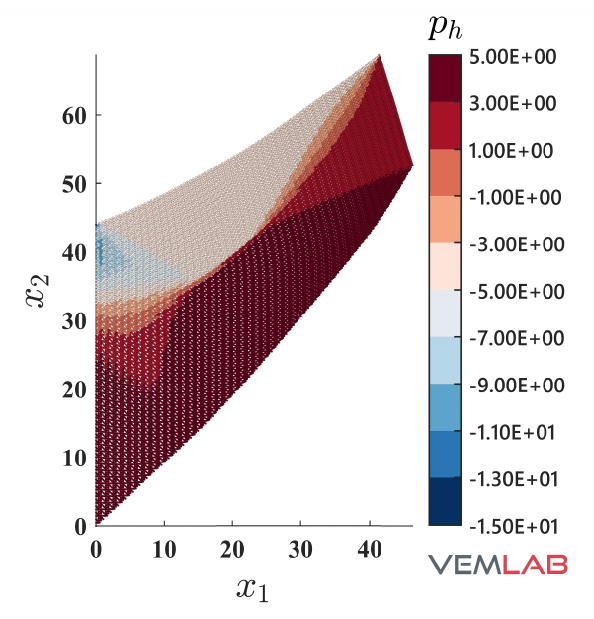}}
\subfigure[] {\label{fig:cook_solutions_b}\includegraphics[width=0.49\linewidth]{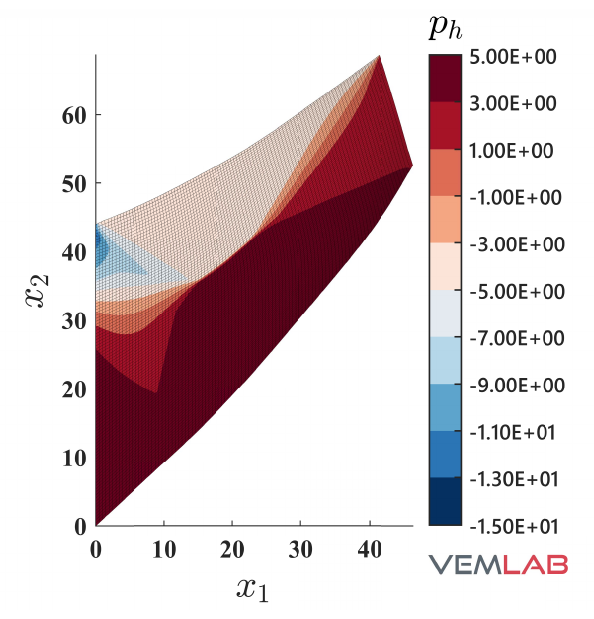}}
}
\mbox{
\subfigure[] {\label{fig:cook_solutions_c}\includegraphics[width=0.49\linewidth]{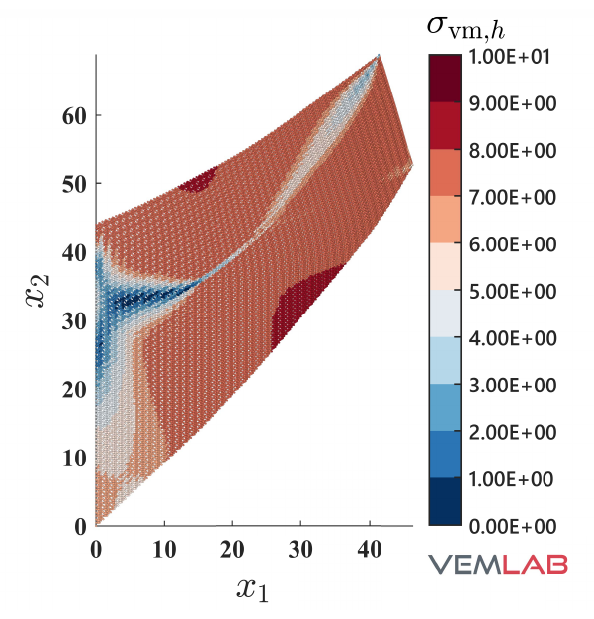}}
\subfigure[] {\label{fig:cook_solutions_c}\includegraphics[width=0.49\linewidth]{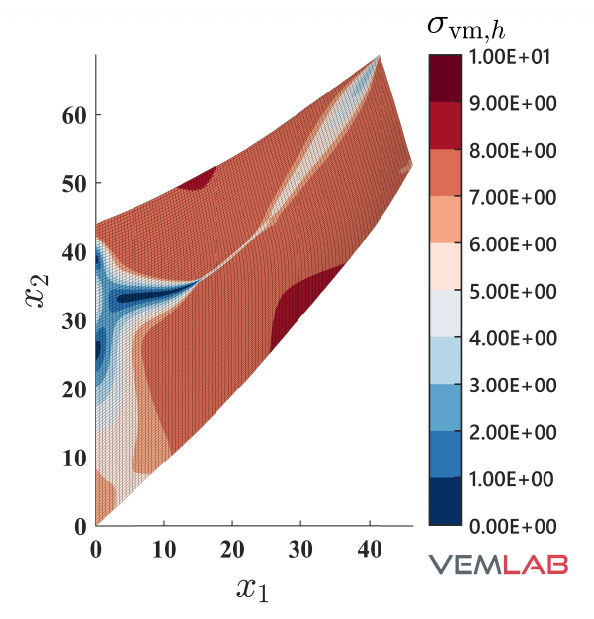}}
}
\caption{Cook's membrane problem. Pressure field solution in MPa ((a) NVEM, (b) FEM Q9 B-bar),
and von Mises stress field solution in MPa ((c) NVEM, (d) FEM Q9 B-bar)}
\label{fig:cook_solutions}
\end{figure}

\fref{fig:cook_solutions} depicts the pressure field and the von Mises 
stress field solutions on the most refined mesh for the NVEM and 
FEM Q9 B-bar approaches. The NVEM plots do not look as smooth as the FEM plots
because of the particular shape of the polygonal element used to construct 
the mesh for the NVEM case (see~\fref{fig:cookmeshes_a}). Despite this, 
the solutions are in good agreement.

\subsection{Tension problem}
\label{sec:tension}

\begin{figure}[!htb]
\centering
\subfigure[] {\label{fig:tensionproblem_a}\includegraphics[width=0.5\linewidth]{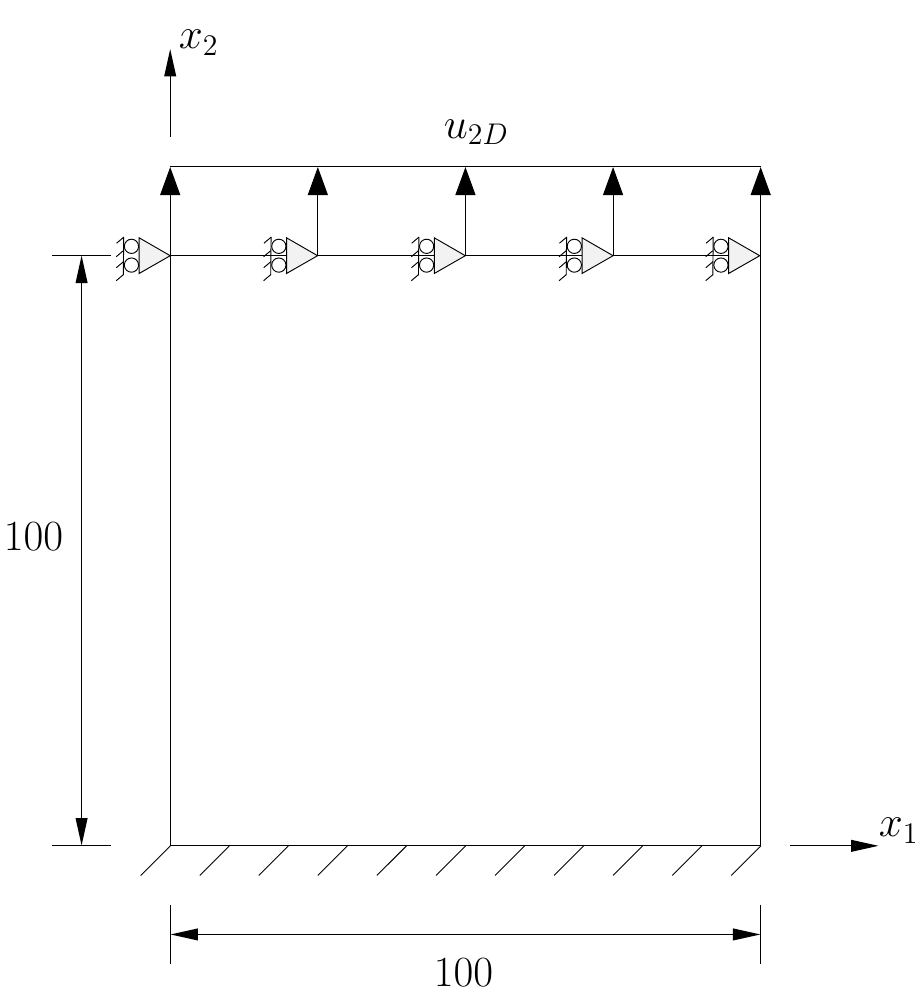}}
\subfigure[] {\label{fig:tensionproblem_b}\includegraphics[width=0.45\linewidth]{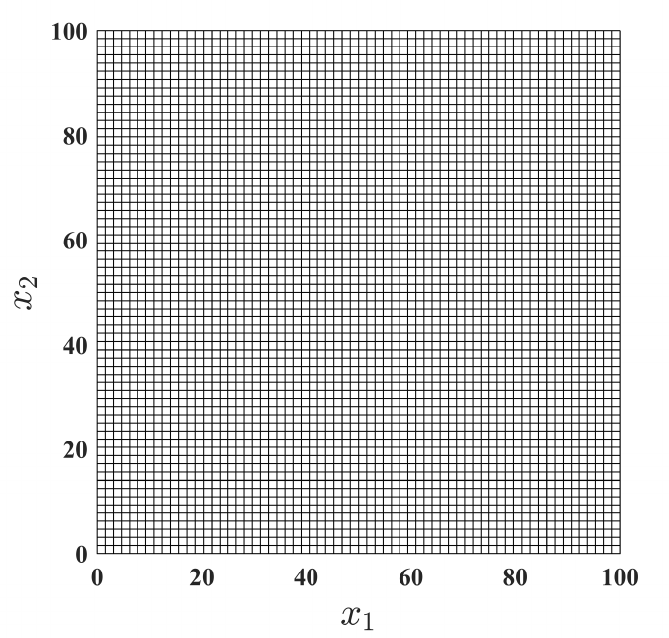}}
\caption{Tension problem. (a) Geometry and boundary conditions, and 
        (b) mesh used for benchmarking the VEM, NVEM, and FEM Q9 B-bar approaches}
\label{fig:tensionproblem}
\end{figure}

In this example, the performance of the NVEM is studied on a pure
tension problem. The domain consists of a square of dimensions 
100 $\times$ 100 mm$^2$ and unit thickness. Plane strain condition is
specified and the following material parameters are used: 
$E_\mathrm{Y}=200000$ MPa, $\nu=0.4999$, $\sigma_{y0}=150$ MPa, $H_i=H_k=0$ MPa 
(perfect plasticity). On the top edge of the domain, a vertical 
displacement of 0.5 mm is imposed while its lateral movement
is restrained. The bottom surface of the domain is fixed. 
\fref{fig:tensionproblem} summarizes the problem definition 
and presents the mesh used in the analysis.

\begin{figure}[!htb]
\centering
\includegraphics[width=0.6\linewidth]{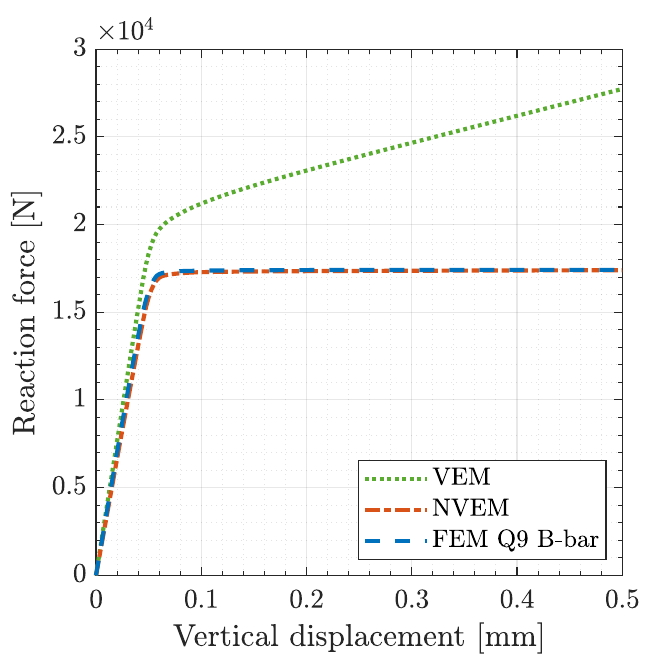}
\caption{Load-displacement curve for the tension problem}
\label{fig:tensionresponsecurve}
\end{figure}

\fref{fig:tensionresponsecurve} presents the load-displacement curve,
where a coincident limit load is observed for the NVEM and FEM Q9 B-bar
approaches. On the other hand, the load-displacement curve for the VEM shows a 
locking effect both in the elastic and plastic regimes.

\begin{figure}[!tbp]
\centering
\mbox{
\subfigure[] {\label{fig:tensionsolutions_a}\includegraphics[width=0.49\linewidth]{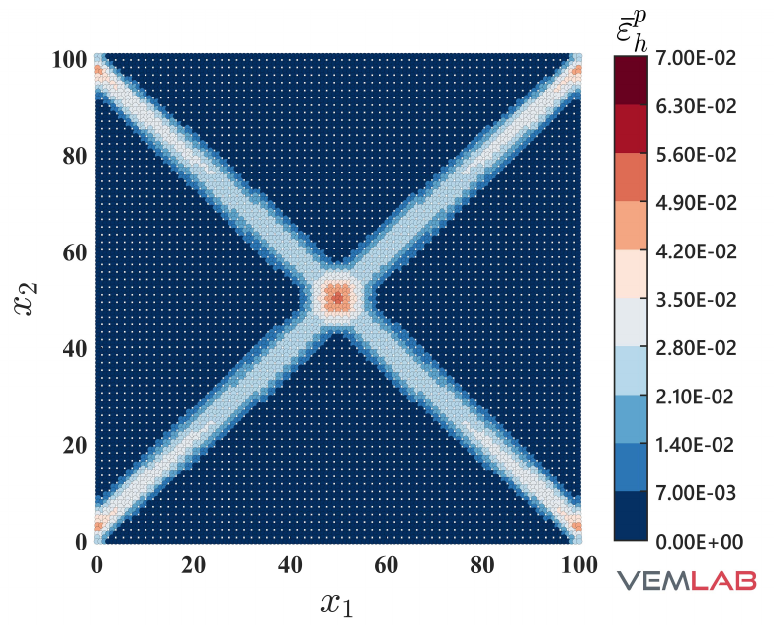}}
\subfigure[] {\label{fig:tensionsolutions_b}\includegraphics[width=0.49\linewidth]{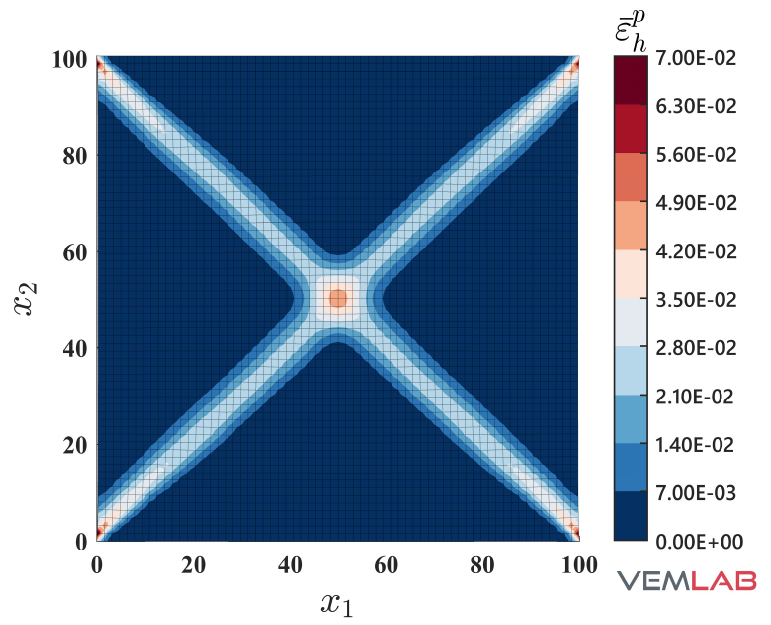}}
}
\mbox{
\subfigure[] {\label{fig:tensionsolutions_c}\includegraphics[width=0.49\linewidth]{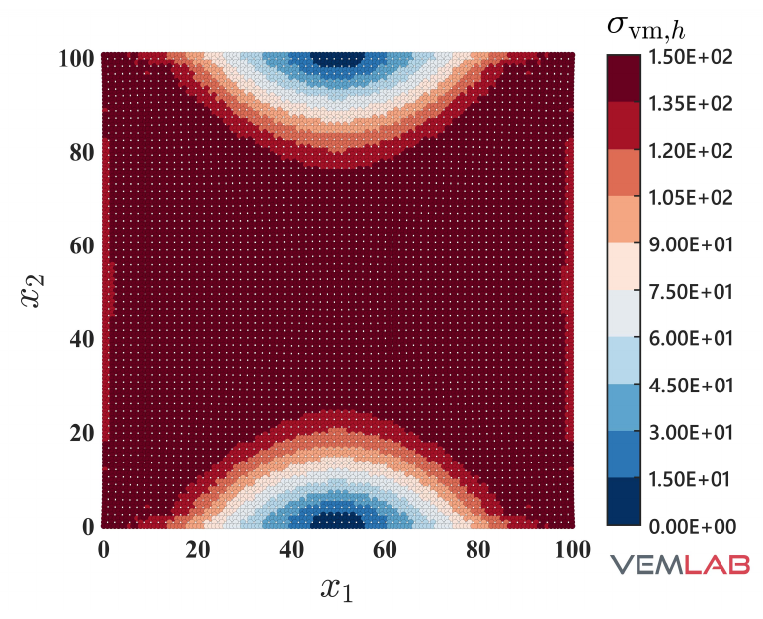}}
\subfigure[] {\label{fig:tensionsolutions_d}\includegraphics[width=0.49\linewidth]{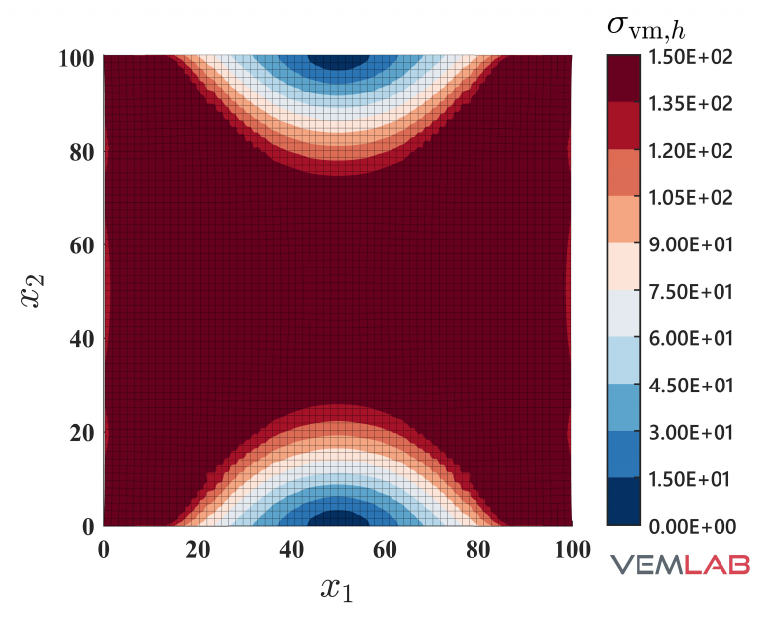}}
}
\caption{Tension problem. Accumulated plastic strain field solution ((a) NVEM, (b) FEM Q9 B-bar),
and von Mises stress field solution in MPa ((c) NVEM, (d) FEM Q9 B-bar)}
\label{fig:tensionsolutions}
\end{figure}

The solutions for the accumulated plastic strain and the von Mises stress
are presented in \fref{fig:tensionsolutions}, where once again good agreement
is observed between the NVEM and FEM Q9 B-bar methods.

\subsection{Perforated plate}
\label{sec:perfplate}

\begin{figure}[!tbp]
\centering
\subfigure[] {\label{fig:perfplatemeshes_a}\includegraphics[width=0.4\linewidth]{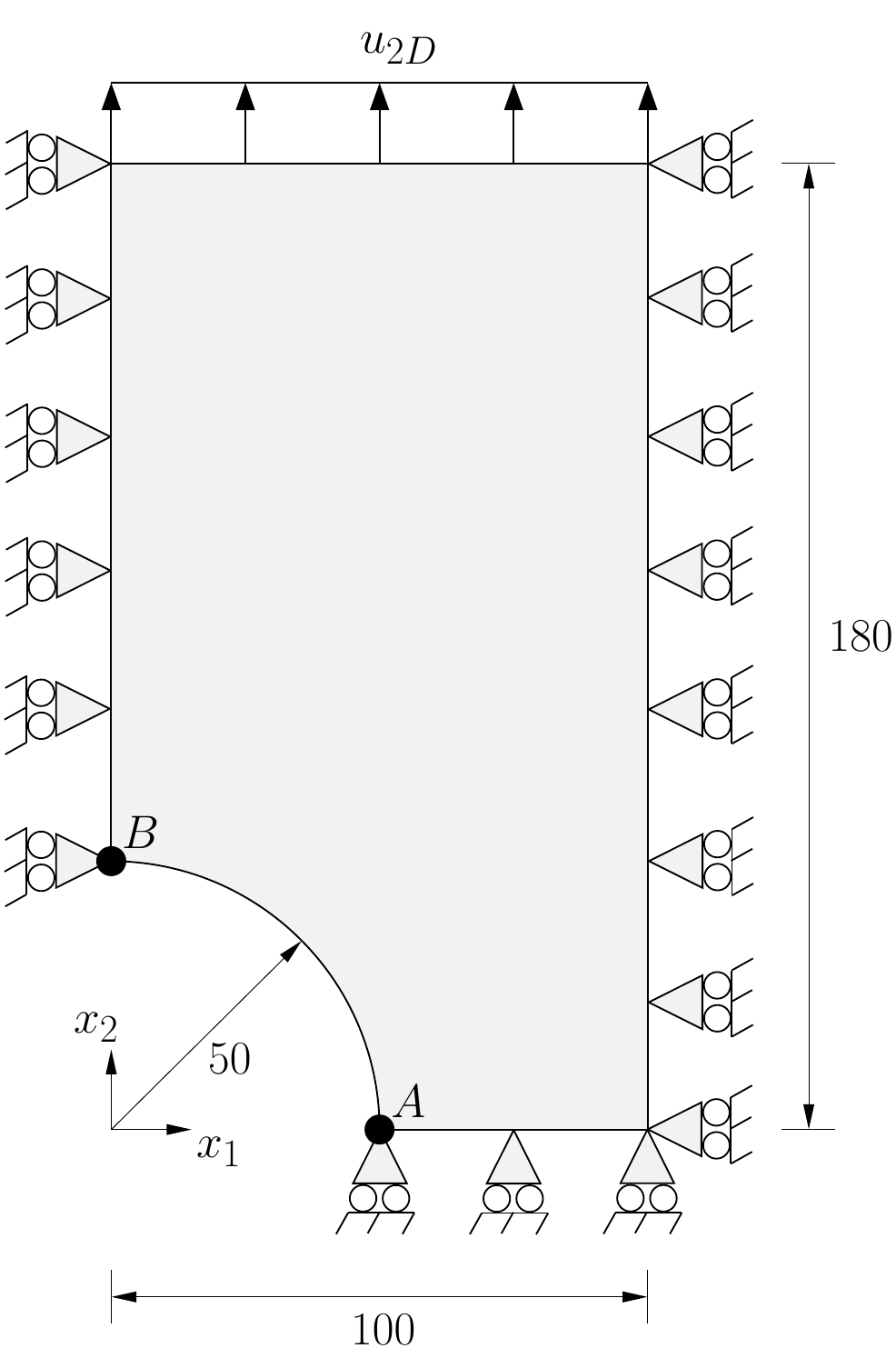}}
\subfigure[] {\label{fig:perfplatemeshes_b}\includegraphics[width=0.37\linewidth]{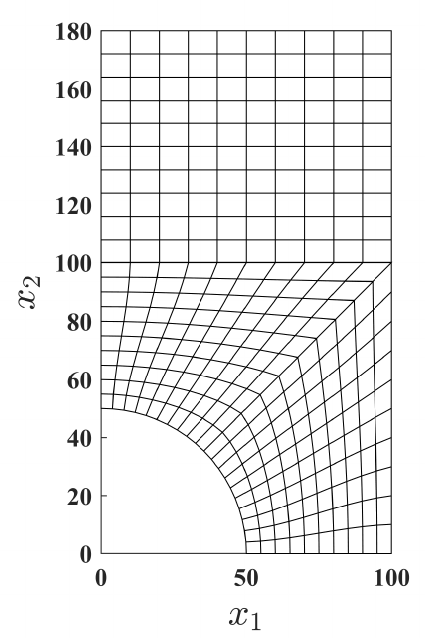}}
\caption{Perforated plate problem. (a) Geometry (with dimensions in mm) and boundary conditions, and 
        (b) a sample mesh used for benchmarking the VEM, NVEM, FEM Q4, and FEM Q9 B-bar approaches}
\label{fig:perfplate_geomeshes}
\end{figure}

\begin{figure}[!htb]
\centering
\includegraphics[width=0.6\linewidth]{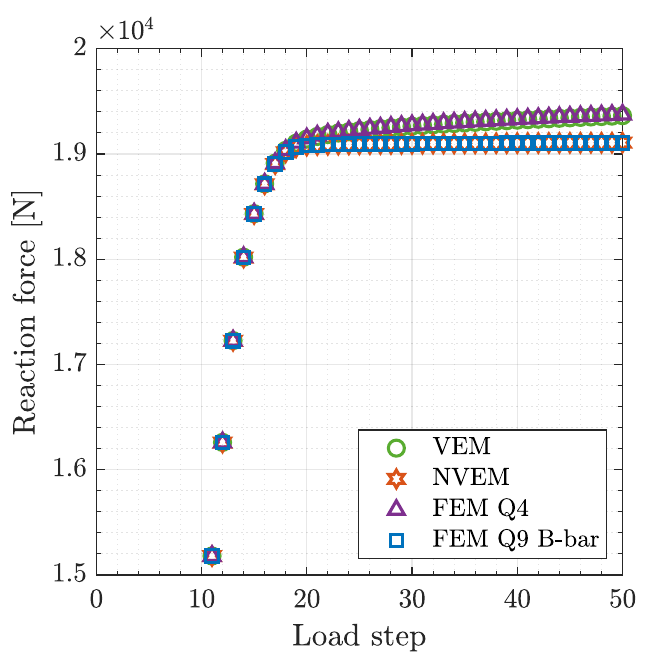}
\caption{Response curve for the perforated plate problem. Reaction force due to the applied
         vertical displacement in steps}
\label{fig:perfplate_responsecurve}
\end{figure}

\begin{figure}[!tbp]
\centering
\mbox{
\subfigure[] {\label{fig:perfplate_solutions_a}\includegraphics[width=0.49\linewidth]{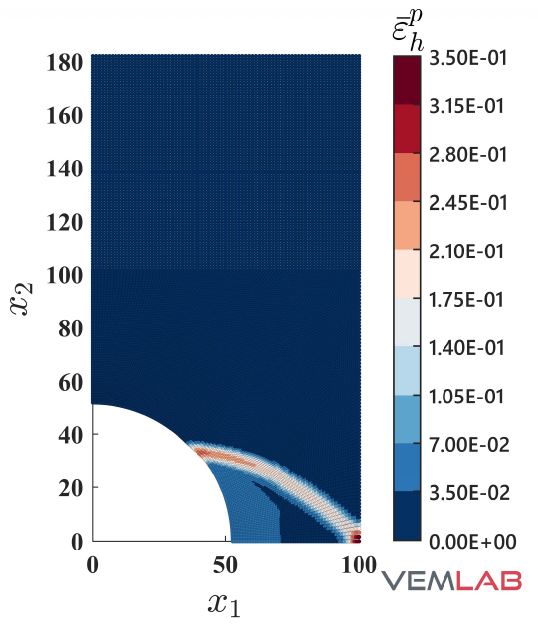}}
\subfigure[] {\label{fig:perfplate_solutions_b}\includegraphics[width=0.49\linewidth]{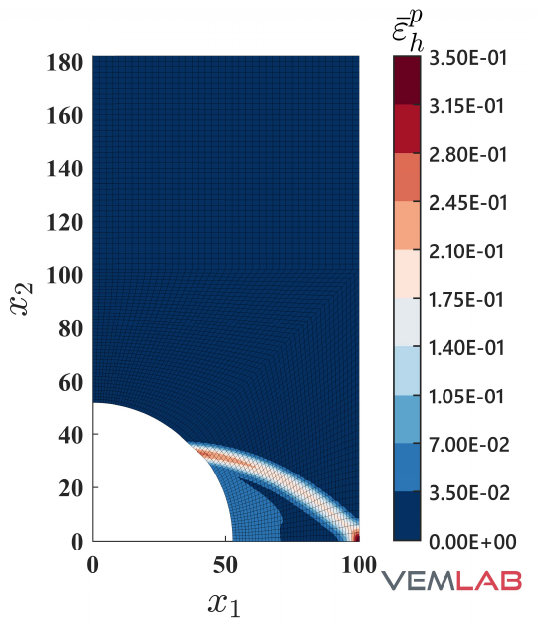}}
}
\mbox{
\subfigure[] {\label{fig:perfplate_solutions_c}\includegraphics[width=0.49\linewidth]{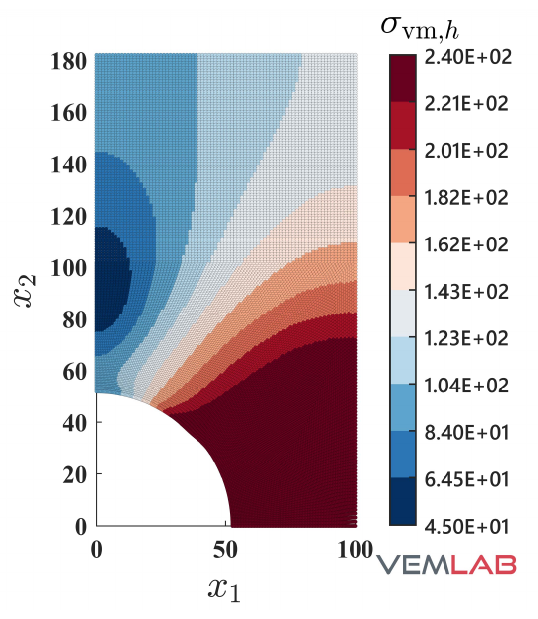}}
\subfigure[] {\label{fig:perfplate_solutions_d}\includegraphics[width=0.49\linewidth]{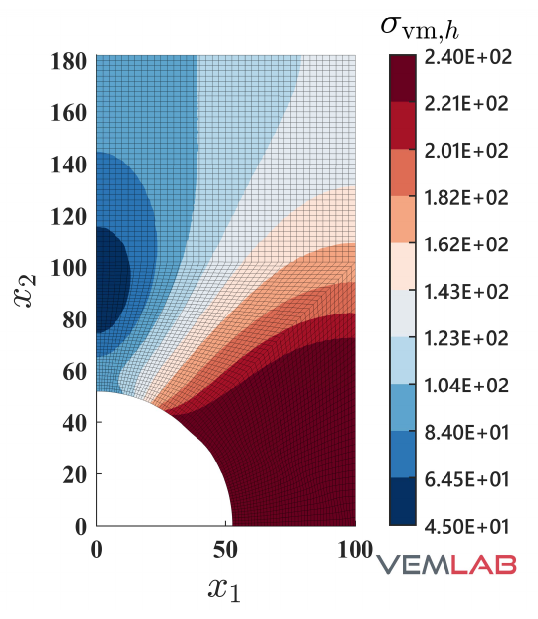}}
}
\caption{Perforated plate problem. Accumulated plastic strain field solution ((a) NVEM, (b) FEM Q9 B-bar),
and von Mises stress field solution in MPa ((c) NVEM, (d) FEM Q9 B-bar)}
\label{fig:perfplate_solutions}
\end{figure}

In this benchmark problem, a quarter of a perforated plate is considered.
The geometry, boundary conditions, and sample mesh are depicted in \fref{fig:perfplate_geomeshes}.
The plate has unit thickness and plane strain condition is assumed.
The material parameters are set to $E_\mathrm{Y}=68646.55$ MPa, $\nu=0.3$, $\sigma_{y0}=238.301595$ MPa,
$H_i=H_k=0$ MPa (perfect plasticity). A vertical displacement $u_{2D}=2$ mm is applied on the top edge of the plate.
On the left, right, and bottom edges of the plate, the translation is restrained in the normal
direction to these edges. Therefore, the plate is highly constrained and thus a locking effect 
is expected for the standard linearly precise VEM. In fact, this is confirmed in 
the response curve shown in \fref{fig:perfplate_responsecurve}, where the expected limit
load is achieved only by the NVEM and the FEM Q9 B-bar in the plastic regime.
A pictorial of the accumulated plastic strain and the von Mises stress 
for the most refined mesh is shown in \fref{fig:perfplate_solutions}, 
where it is observed that the NVEM and FEM Q9 B-bar solutions look very similar.

A comparison of the horizontal displacement at point $A$ and the
vertical displacement at point $B$ among the different methods
is presented in Table \ref{table:perfplate}, where a perfect match
to three decimal places is obtained for the NVEM and FEM Q9 B-bar methods.

\begin{table}[htbp]
  \centering
  \captionbox{Displacement comparison for the perforated plate problem\label{table:perfplate}}{
    \begin{tabular}{lll}
    \midrule
    \textbf{Method} & $\bm{u_{1A}}$ \textbf{[mm]} & $\bm{u_{2B}}$ \textbf{[mm]}\\
    \midrule
    VEM   & 2.716 & 1.850 \\
    NVEM  & 2.745 & 1.855 \\
    FEM Q4 & 2.714 & 1.849 \\
    FEM Q9 B-bar & 2.745 & 1.855 \\
    \midrule
    \end{tabular}}
\end{table}

\subsection{Prandtl's punch test}
\label{sec:prandtl}

\begin{figure}[!tbp]
\centering
\includegraphics[width=0.9\linewidth]{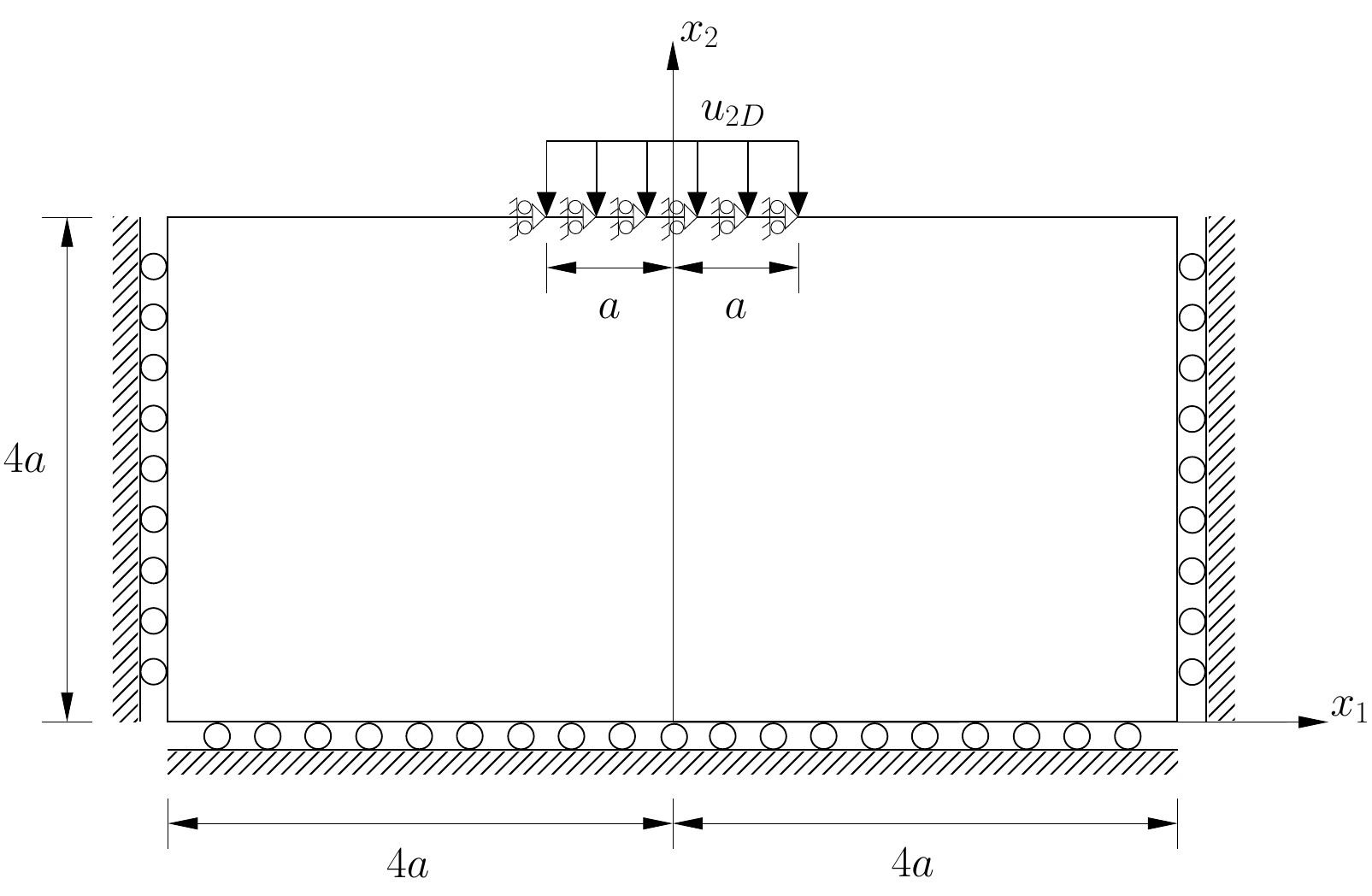}
\caption{Geometry and boundary conditions for the Prandtl's punch test}
\label{fig:prandtlproblem}
\end{figure}

\begin{figure}[!tbp]
\centering
\mbox{
\subfigure[] {\label{fig:prandtlmeshes_a}\includegraphics[width=0.49\linewidth]{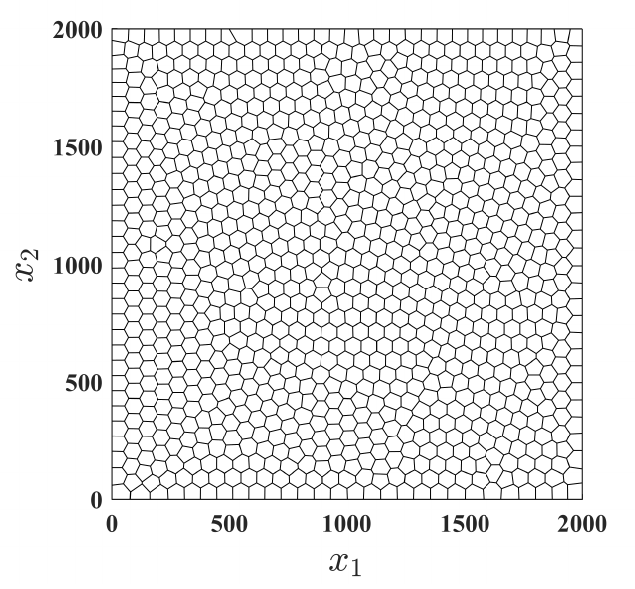}}
\subfigure[] {\label{fig:prandtlmeshes_b}\includegraphics[width=0.49\linewidth]{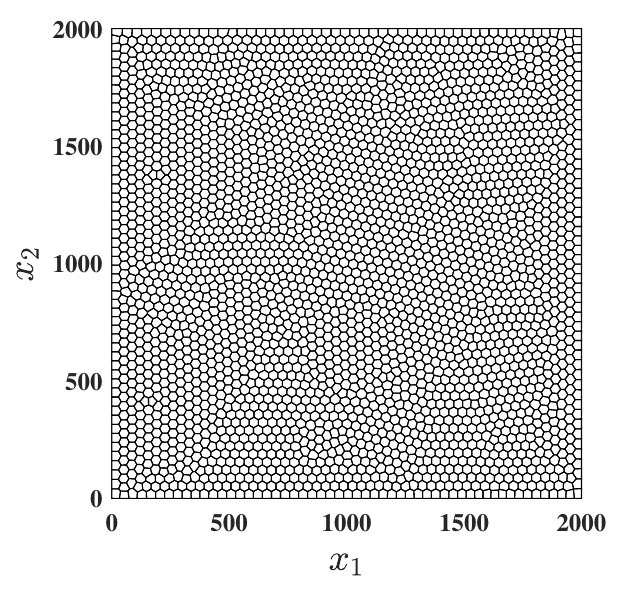}}
}
\mbox{
\subfigure[] {\label{fig:prandtlmeshes_c}\includegraphics[width=0.49\linewidth]{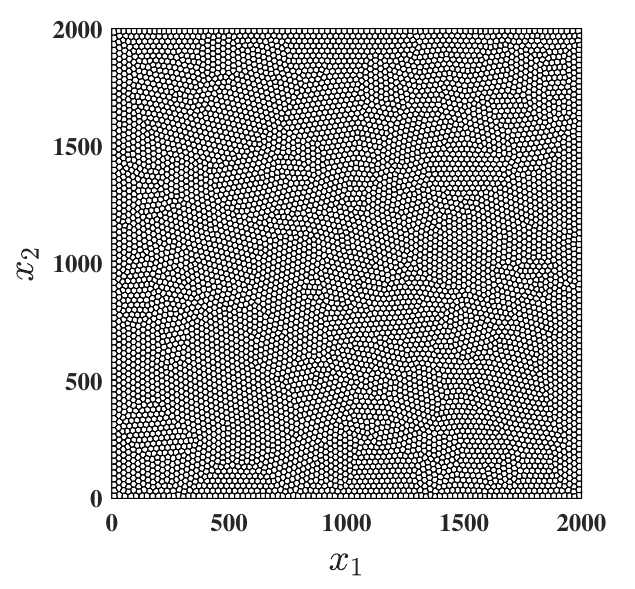}}
\subfigure[] {\label{fig:prandtlmeshes_d}\includegraphics[width=0.49\linewidth]{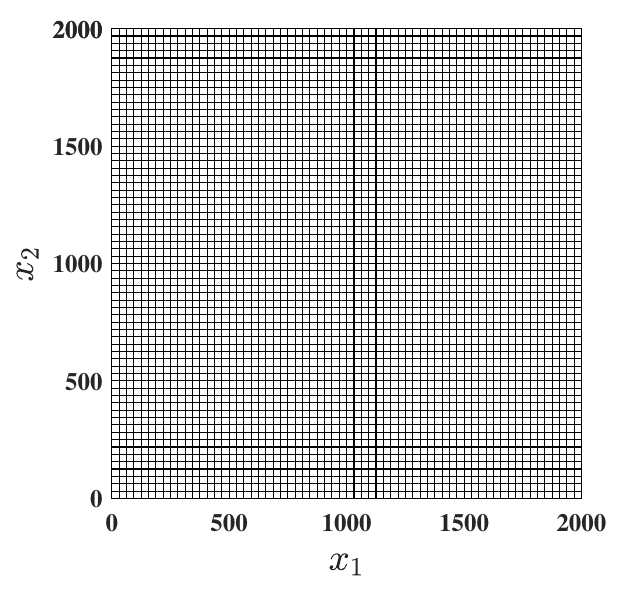}}        
}
\caption{Meshes for the Prandtl's punch test. (a) 4004 DOF polygonal mesh, (b) 12004 DOF polygonal mesh, 
        (c) 32004 DOF polygonal mesh, and (d) 33282 DOF 9-node quadrilateral mesh. Meshes (a)-(c) are used
        for the VEM and NVEM approaches, whereas mesh (d) is used for the FEM Q9 B-bar approach}
\label{fig:prandtlmeshes}
\end{figure}

\begin{figure}[!htb]
\centering
\includegraphics[width=0.6\linewidth]{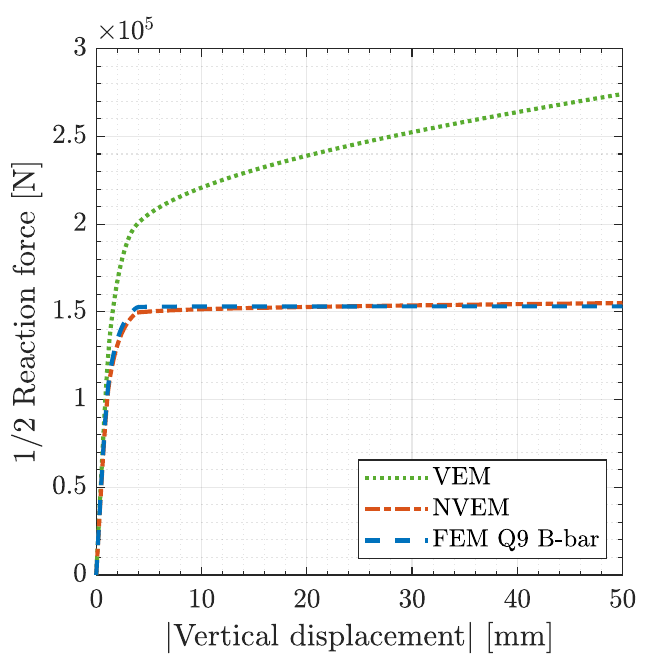}
\caption{Load-displacement curve for the Prandtl's punch test}
\label{fig:prandtlresponsecurve}
\end{figure}

\begin{figure}[!tbp]
\centering
\mbox{
\subfigure[] {\label{fig:prandtlpstrains_a}\includegraphics[width=0.49\linewidth]{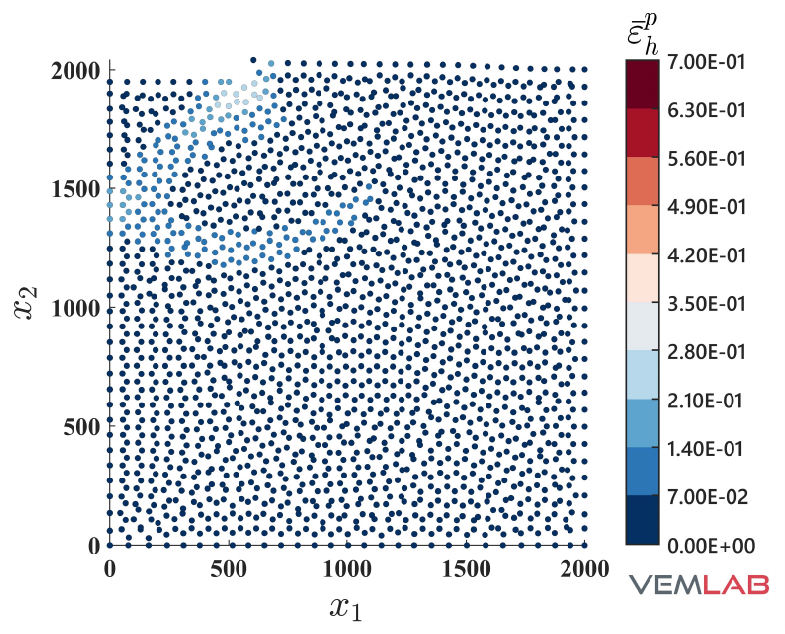}}
\subfigure[] {\label{fig:prandtlpstrains_b}\includegraphics[width=0.49\linewidth]{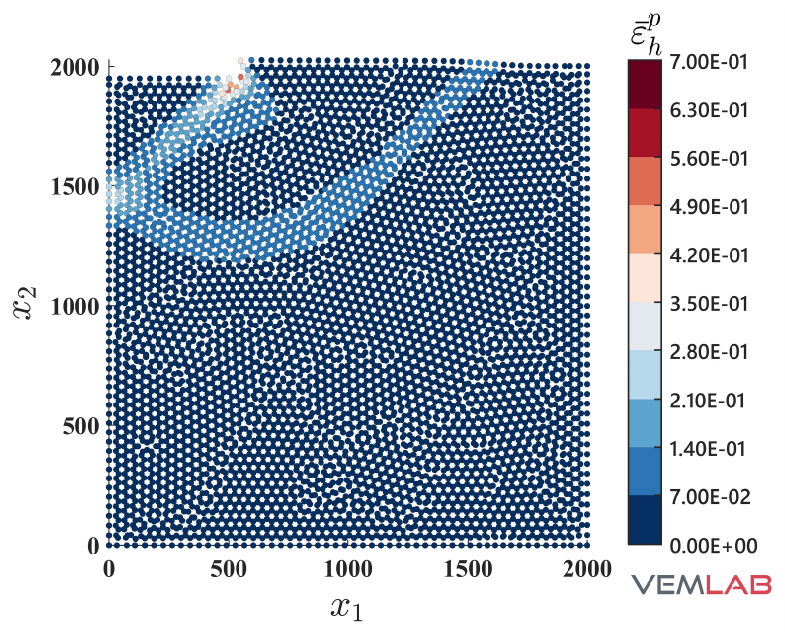}}
}
\mbox{
\subfigure[] {\label{fig:prandtlpstrains_c}\includegraphics[width=0.49\linewidth]{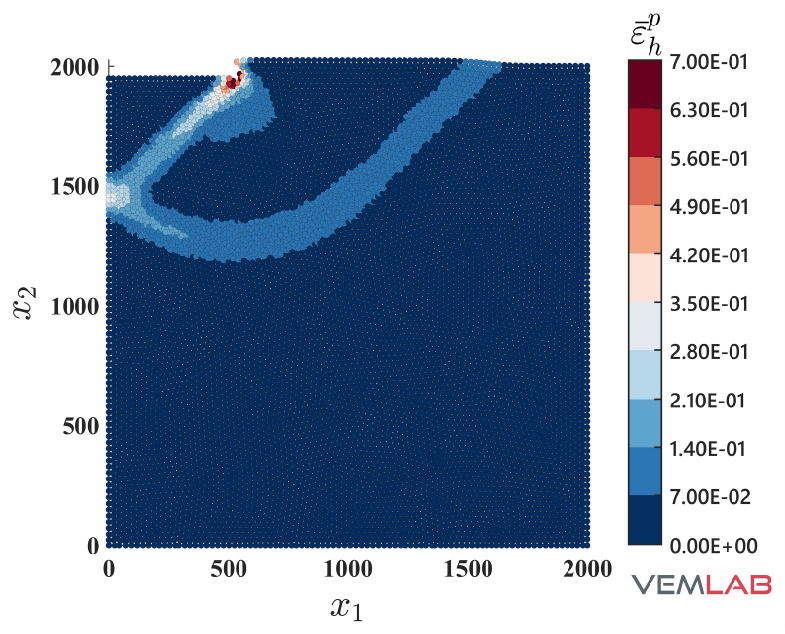}}
\subfigure[] {\label{fig:prandtlpstrains_d}\includegraphics[width=0.49\linewidth]{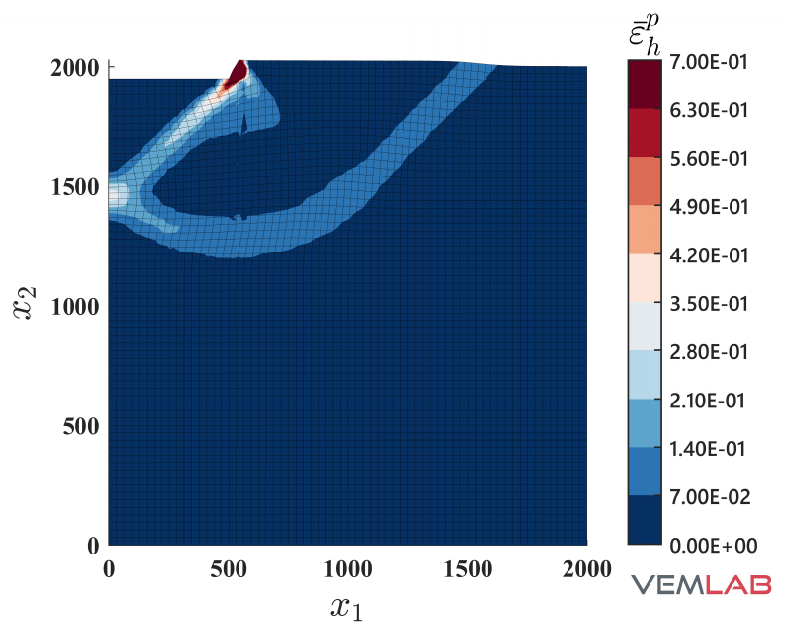}}
}
\caption{Accumulated plastic strain field solution for the Prandtl's punch test. (a) NVEM (4004 DOF polygonal mesh), (b) NVEM (12004 DOF polygonal mesh), 
        (c) NVEM (32004 DOF polygonal mesh), and (d) FEM Q9 B-bar (33282 DOF 9-node quadrilateral mesh)}
\label{fig:prandtlpstrains}
\end{figure}

\begin{figure}[!tbp]
\centering
\mbox{
\subfigure[] {\label{fig:prandtlpressure_a}\includegraphics[width=0.49\linewidth]{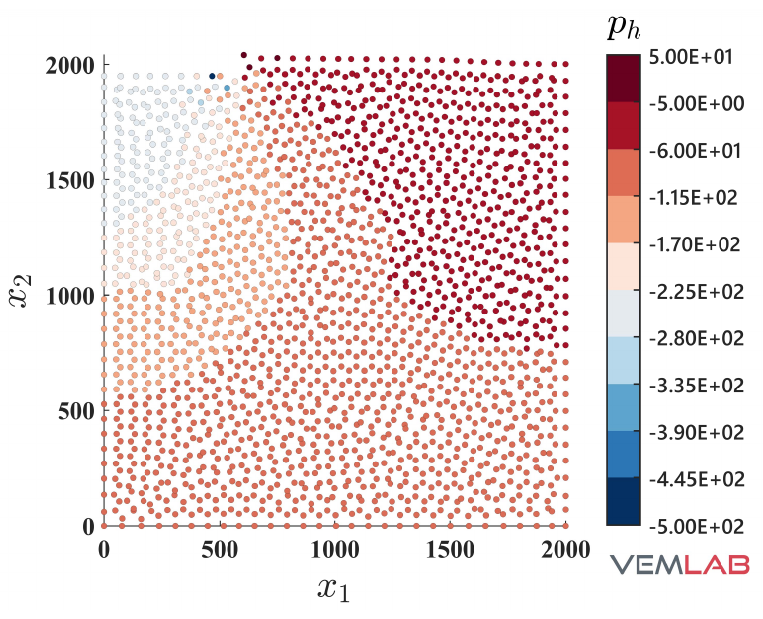}}
\subfigure[] {\label{fig:prandtlpressure_b}\includegraphics[width=0.49\linewidth]{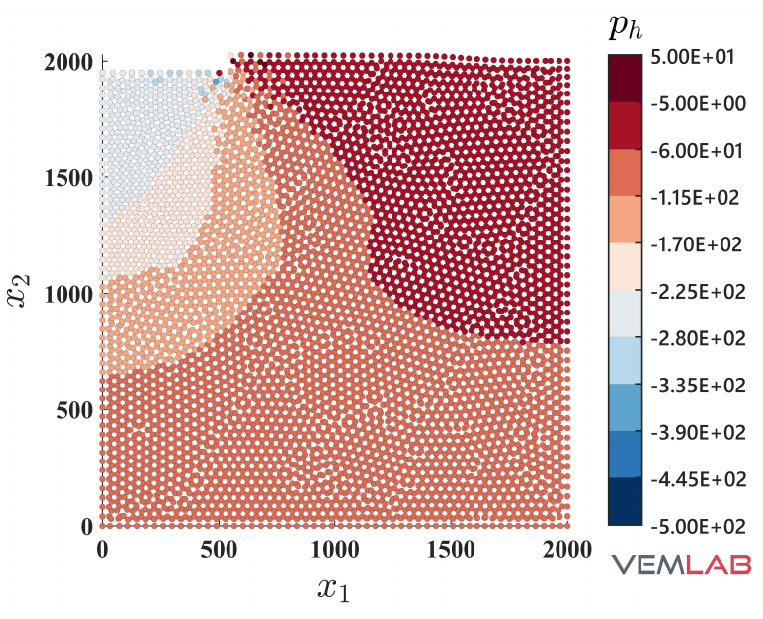}}
}
\mbox{
\subfigure[] {\label{fig:prandtlpressure_c}\includegraphics[width=0.49\linewidth]{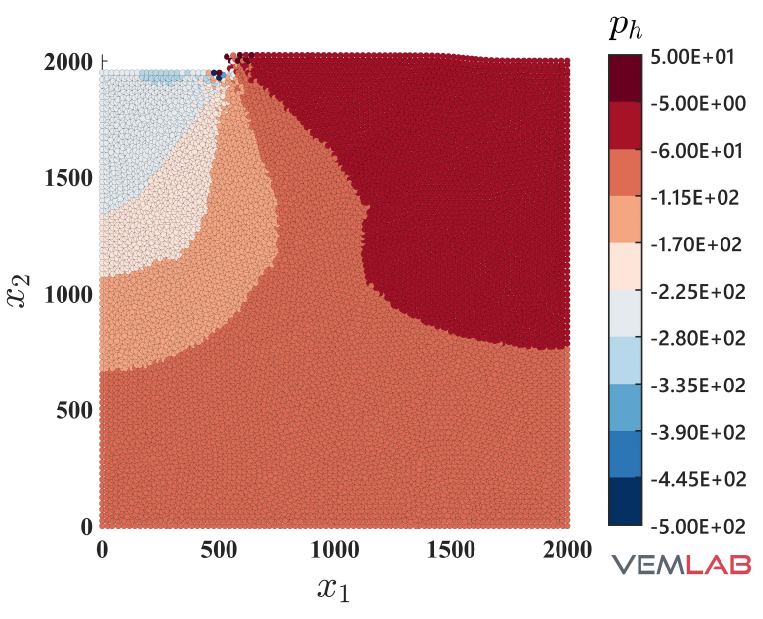}}
\subfigure[] {\label{fig:prandtlpressure_d}\includegraphics[width=0.49\linewidth]{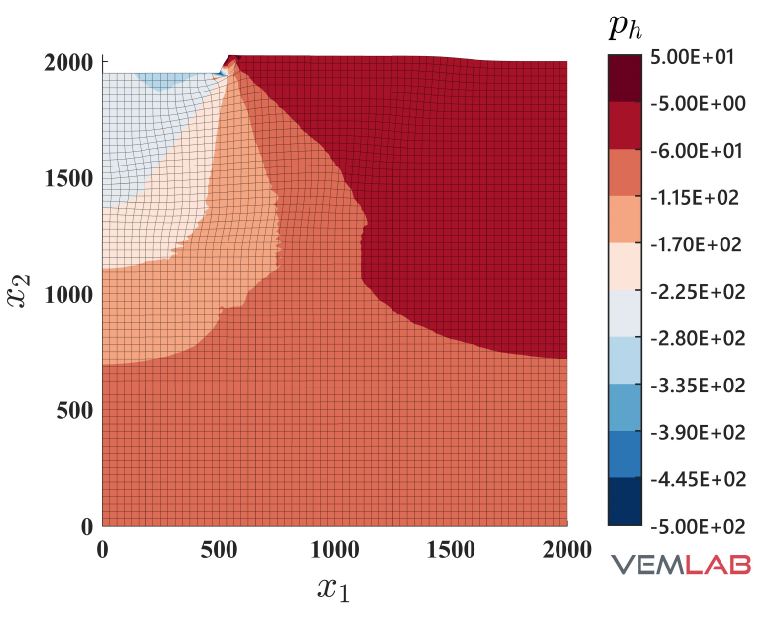}}
}
\caption{Pressure field solution in MPa for the Prandtl's punch test. (a) NVEM (4004 DOF polygonal mesh), (b) NVEM (12004 DOF polygonal mesh), 
        (c) NVEM (32004 DOF polygonal mesh), and (d) FEM Q9 B-bar (33282 DOF 9-node quadrilateral mesh)}
\label{fig:prandtlpressure}
\end{figure}

\begin{figure}[!tbp]
\centering
\subfigure[] {\label{fig:prandtlvemresults_a}\includegraphics[width=0.49\linewidth]{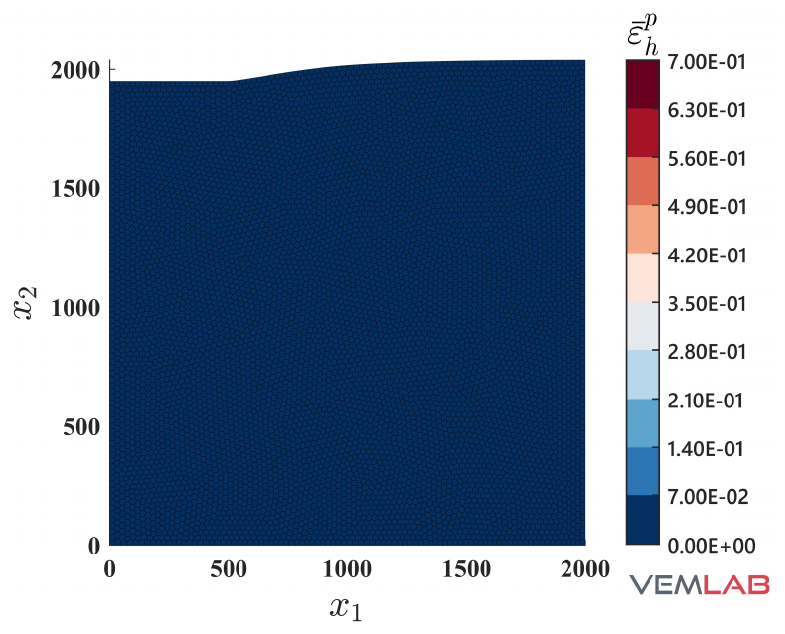}}
\subfigure[] {\label{fig:prandtlvemresults_b}\includegraphics[width=0.49\linewidth]{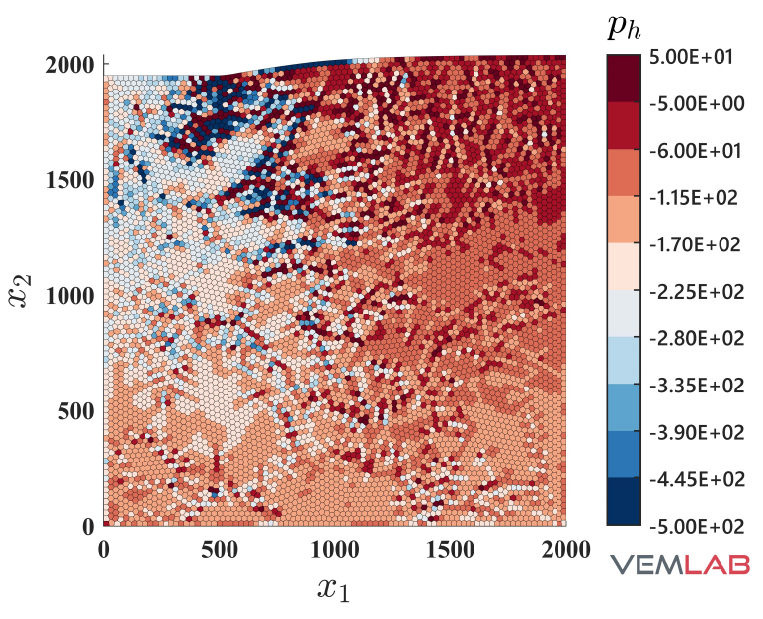}}
\caption{Prandtl's punch test. VEM solution (32004 DOF polygonal mesh) for the (a) accumulated plastic strain, and (b) pressure field in MPa}
\label{fig:prandtlvemresults}
\end{figure}

The last benchmark problem is devoted to demonstrate the performance of
the NVEM in a highly constrained compression problem. The Prandtl's
punch test~\cite{Tan-Li-Chou:1989} is selected for this purpose. The geometry and
boundary conditions are shown in~\fref{fig:prandtlproblem}, where
the punch surface located on the top is horizontally restrained (rough punch)
while a downward vertical displacement $u_{2D}=-50$ mm is imposed on it.
The dimensions are defined using $a=500$ mm. Because of the symmetry, only
half of the domain is considered for discretization. Three polygonal
meshes with increasing refinements are considered for the VEM and 
NVEM (\fref{fig:prandtlmeshes}(a)--(c)). The most refined VEM/NVEM mesh 
(\fref{fig:prandtlmeshes_c}) and the mesh for the FEM Q9 B-bar (\fref{fig:prandtlmeshes_d})
have similar number of DOF. Unit thickness is considered and 
plane strain condition is assumed with material parameters 
set to $E_\mathrm{Y}=10^5$ MPa, $\nu=0.499$, $\sigma_{y0}=E_\mathrm{Y}/1000$ MPa,
$H_i=H_k=0$ MPa (perfect plasticity). As observed in \fref{fig:prandtlresponsecurve}, 
a severe locking behavior is obtained for the standard linearly precise VEM
in this highly constrained problem. On the other hand, the same figure
reports the nearly coincident limit load that is obtained for the NVEM and FEM Q9 B-bar
approaches. The accumulated plastic strain is reported on \fref{fig:prandtlpstrains},
where similar results are obtained for the NVEM and FEM Q9 B-bar on the
meshes with similar number of DOF (\fref{fig:prandtlpstrains_c}
and \fref{fig:prandtlpstrains_d}, respectively). The pressure
field is depicted in \fref{fig:prandtlpressure}, where the NVEM and FEM Q9 B-bar
solutions look very similar and smooth on the meshes with similar number of 
DOF (\fref{fig:prandtlpressure_c} and \fref{fig:prandtlpressure_d}, respectively). 
Just for completeness, the accumulated plastic strain and the pressure field solutions for the 
linearly precise VEM are shown in \fref{fig:prandtlvemresults}, where 
the severe locking behavior is clearly observed in the oscillations of
the pressure field.

\section{Summary and conclusions}
\label{sec:conclusions}

The node-based uniform strain virtual element method (NVEM) that was recently
proposed for compressible and nearly incompressible elasticity~\cite{Ortiz-Silva-Salinas-Hitschfeld-Luza-Rebolledo:2023} 
has been extended to elastoplastic solids at small strains.
In the proposed method, the strain is averaged at the nodes from the 
strain of surrounding linearly precise virtual elements using 
a generalization to virtual elements of the node-based uniform strain 
approach for finite elements~\cite{Dohrmann-Heinstein-Jung-Key-Witkowski:2000}. 
The averaged strain is then used to sample the weak form at the nodes 
of the mesh leading to a method in which all the field variables, including state 
and history-dependent variables, are related to the nodes. Consequently, 
in the nonlinear computations these variables are tracked only at the nodal locations.

Various elastoplastic benchmark problems were conducted to assess 
the performance of the NVEM. These included a thick-walled cylinder
under internal pressure, a combined bending and shear problem (Cook's membrane), 
a pure tension problem, a perforated plate subjected to a displacement
producing tension, and a highly constrained problem in compression (Prandtl's 
punch test). The comparisons with the well-known locking-free 
9-node B-bar quadrilateral finite element~\cite{Simo-Hughes:1998} 
revealed that the NVEM effectively enables linearly precise 
virtual elements to solve elastoplastic solids with 
accuracy and is locking-free. For perfect plasticity, these comparisons
also demonstrated that the NVEM is able to capture the expected limit load.
Finally, we mention that the present work completes our short term scope 
for the NVEM development and that its extension to large deformations 
with remeshing is an undergoing work.

\backmatter

\bmhead{Acknowledgements}

This work was performed under the auspices of the Chilean National Fund for Scientific and Technological Development
(FONDECYT) through grant ANID FONDECYT No. 1221325 (R.S-V and A.O-B).









\clearpage


\end{document}